\documentclass[a4paper,12pt,twoside]{elsarticle}
\pdfoutput=1
\usepackage[utf8]{inputenc}
\usepackage{amssymb,amsmath,amsthm}
\usepackage{amscd}
\usepackage{a4wide}
\usepackage{setspace}
\usepackage{multirow}
\usepackage{color}
\usepackage{natbib}
\usepackage{pdfpages}
\usepackage{booktabs}
\usepackage{graphics}
\usepackage{subfigure}

\newcommand{\R}{\mathbb{R}}

\theoremstyle{definition}

\theoremstyle{definition}

\theoremstyle{plain}
\newtheorem{proposition}{Proposition}[section]

\newtheorem{theorem}{Theorem}

\newcommand{\bn}{\mathbf{b}}
\newcommand{\xn}{\mathbf{x}}
\newcommand{\un}{\mathbf{u}}
\newcommand{\vn}{\mathbf{v}}
\newcommand{\rn}{\mathbf{r}}
\newcommand{\yn}{\mathbf{y}}
\newcommand{\pn}{\mathbf{p}}

\newcommand{\fn}{\mathbf{f}}

\newcommand{\gn}{\mathbf{g}}

\newcommand{\nn}{\mathbf{n}}
\newcommand{\wn}{\mathbf{w}}
\newcommand{\An}{\mathbf{A}}

\newcommand{\Fn}{\mathbf{F}}

\newcommand{\Hn}{\mathbf{H}}
\newcommand{\Ln}{\mathbf{L}}

\begin{document}

\begin{frontmatter}



\title{Performance Evaluation of Block Diagonal Preconditioners for Divergence-Conforming B-Spline Discretization of Stokes Equations}


\author[1,2]{Adriano M. A. C\^{o}rtes} 
\author[1]{Alvaro L. G. A. Coutinho \corref{cor1}}
\cortext[cor1]{Correponcence to: A.L.G.A Coutinho, High Performance Computing Center and Department of Civil Engineering, COPPE/Federal University of Rio de Janeiro, PO Box 68506, Rio de Janeiro, RJ 21945, Brazil.}
\ead{alvaro@nacad.ufrj.br}

\address[1]{High Performance Computing Center and Department of Civil Engineering, COPPE, Federal University of Rio de Janeiro, Rio de Janeiro, Brazil}
\address[2]{Department of Mathematics and Statistics, Federal State University of Rio de Janeiro, Rio de Janeiro, Brazil}

\begin{abstract}
The recently introduced divergence-conforming B-spline discretizations allow the construction of smooth discrete velocity-pressure pairs for viscous incompressible flows that are at the same time inf-sup stable and divergence-free. When applied to discretize Stokes equations, it generates a symmetric and indefinite linear system of saddle-point type. Krylov subspace methods are usually the most efficient procedures to solve such systems. One of such methods, for symmetric systems, is the Minimum Residual Method (MINRES). However, the efficiency and robustness of Krylov subspace methods is closely tied to appropriate preconditioning strategies. For the discrete Stokes system, in particular, block-diagonal strategies provide efficient preconditioners. In this paper, we compare the performance of block-diagonal preconditioners for several block choices. We verify how eigenvalue clustering promoted by the preconditioning strategies affects MINRES convergence. We also compare the number of iterations and wall-clock timings. We conclude that an incomplete Cholesky block-diagonal preconditioning strategy with relaxed inner conjugate gradients iterations provides the best computational strategy when compared to other block-diagonal and global solution strategies.
\end{abstract}


\end{frontmatter}

\section{Introduction}

The concept of Isogeometric Analysis (IGA) first appeared in \cite{HBC2005}, and since then several papers followed, either exploring their mathematical theory, for example, \cite{Bazilevs_Est} and \cite{Sangalli_Est}, or showing their potential in engineering applications, e.g. \cite{Cottrell2006}, \cite{Bazilevs2006}, \cite{Bazilevs2007}, \cite{Zhang2007}, \cite{Gomez2008}, \cite{Bazilevs2008}, \cite{Benson2010}, \cite{Bazilevs2010} and \cite{Akkerman2011}. In \cite{BuffaSV2010}, the IGA concept was used to discretize vector fields of electro-magnetics problems. For such problems,  it is known that the function spaces satisfy a de Rham diagram at continuous level, and, for a discretization to be successfully applied to them, the finite dimensional spaces should also satisfy the de Rham diagram at the discrete level. Exploring one of the main features of spline basis functions, that is the easy control of the basis polynomial degree and regularity, and by a suitable choice of B-spline spaces of each component of the two-dimensional vector field, Buffa et al. \cite{BuffaSV2010} introduced an IGA discretization satisfying a de Rham diagram. They have shown that the technique can be viewed as a smooth generalization of Nédélec elements, and thus good results were reported.

The generalization for three-dimensional vector fields and the mathematical theory of such discretization appeared in \cite{BuffaRSV2011}. Their approach, called Isogeometric Discrete Differential Forms, was inspired by the theory of finite element exterior calculus of Arnold et al. \cite{Arnold2010}. The main ingredient of such discretization scheme is the devising of suitable projectors to the continuous spaces onto the discrete ones, that render commuting diagrams and the satisfaction of a de Rham diagram at the discrete level. As a consequence, this implies (numerical) stability.

Three similar vector discretizations for the Stokes problem were introduced in \cite{BuffaFS2011}. By a proper choice of the polynomial degrees and the regularity of the discrete velocity field components and the discrete pressure field, these discretizations can be interpreted as smooth generalizations of Nédélec, Taylor-Hood and Raviart-Thomas elements. Because of the smoothness of the basis functions, the discrete velocity spaces of these elements are $\mathbf{H}^1$-conforming, which make them suitable to discretize the Stokes system. Furthermore, in the case of the Raviart-Thomas element type, Buffa et al. \cite{BuffaFS2011} characterize the image of the divergence operator from the discrete velocity space (with and without boundary conditions) onto the discrete pressure space. In this way, by choosing the discrete velocity-pressure pair such that the divergence operator is surjective, this element provides a point-wise divergence-free discrete vector field, a condition that is generally only satisfied weakly by classical mixed finite elements.

Following the developments of Buffa et al. \cite{BuffaRSV2011} and also inspired by the theory of finite element exterior calculus of Arnold et al. \cite{Arnold2010}, Evans \cite{Evans2011} further developed the Raviart-Thomas element type in the context of Hilbert complexes. Indeed, by using the stable projectors of \cite{BuffaRSV2011}, a divergence-preserving transformation (Piola transformation) of the velocity field and an integral-preserving transformation of the pressure field, Evans devised a Stokes complex with a compatible sub-complex that furnishes a discretization scheme, that is, at the same time {\it inf-sup} stable and divergence-free. In \cite{Evans2011}, Evans applied this discretization scheme for several viscous incompressible flows and also developed their mathematical theory, as well. We will overview this discretization scheme and part of its mathematical theory in section \ref{sec:divconf}.

It is well known that the discretization of the Stokes equations by {\it inf-sup} stable mixed elements requires the solution of a symmetric indefinite linear system, called the (discrete) Stokes system, with a block coefficient matrix of saddle-point type. Several strategies for solving the Stokes linear system have appeared in the literature (\cite{Elman94multigridand}, \cite{Peters_fastiterative}, \cite{elman2005finite} and \cite{Benzi2005survey}), the most popular being variants of Uzawa method, like {\it inexact} Uzawa method, and the Minimum Residual Method (MINRES) \cite{Paige1975}. The latter is a member from the Krylov subspace methods family, and as such, its robustness and performance is highly dependent of a good preconditioning strategy. Among the possibilities for the Stokes system are the block-diagonal preconditioners introduced by Wathen and Silvester \cite{Wathen1993}, \cite{Wathen1994}.

It is necessary to mention that MINRES is being used to solve large-scale problems in science, such as earth mantle convection flows in parallel by finite elements with octree-based adaptive mesh refinement and coarsening (AMR/C), demonstrating scalability up to 122880 cores \cite{Burstedde2013}. For the same application, in \cite{Kronbichler2012} is shown that efficient block solvers, based on the combination of flexible GMRES and algebraic multigrid, present an almost invariant number of iterations with mesh refinement.


The remainder of this paper is organized as follows. In section 2, we review some isogeometric analysis definitions, in order to setup the nomenclature for the divergence-conforming discretization. Section 3 is dedicated to reviewing the results of \cite{Evans2011} and \cite{Evans2012} with respect to Stokes flow. First we present the discrete velocity-pressure pair on the parametric domain, and how it is mapped to general geometries by means of proper transformations. Also, the inf-sup stability and the divergence-free property of the divergence-conforming discrete velocity-pressure pair is presented. The next section deals with the discrete variational problem, and how Nitsche's method is used to impose Dirichlet boundary condition weakly. In section 5, we review the Minimum Residual Method, mainly its convergence properties and related preconditioners. We also present the block-diagonal preconditioning strategy of Wathen and Silvester \cite{Wathen1993}, \cite{Wathen1994}, and the choices we made for the blocks. Section 6 is devoted to numerical results. We present the results for three examples: two manufactured analytical solutions, for different geometries and the lid-driven cavity flow benchmark. For the last one, we do a detailed analysis of the preconditioners performance. The paper ends with a summary of main conclusions.

\section{Isogeometric standards: Spline spaces and the geometrical mapping}
With the purpose of defining the divergence-conforming discretization of \cite{Evans2011}, in the next section, we need to recall some spline space definitions and notations. Here, we follow closely \cite{BuffaFS2011}, \cite{BuffaRSV2011}, \cite{Evans2011} and \cite{Evans2012}.

\subsection{Univariate B-Splines}
To define a univariate B-spline basis one needs to specify the number $n$ of basis functions wanted, the polynomial degree $p$ of the basis and a knot vector $\Xi$. A knot vector $\Xi$ is a finite nondecreasing sequence $\Xi = \{ 0=\xi_1, \ldots, \xi_{n+p+1}=1 \}$. They can have repeated knots, in this case one says that the knot has multiplicity greater than one. Introducing the vector $\boldsymbol{\zeta} = \{\zeta_1,\ldots,\zeta_m\}$ of knots without repetitions, also called breakpoints, and the vector $\{r_1,\ldots,r_m\}$ of their corresponding multiplicities, one has that,
\begin{equation}
	\Xi = \{ \underbrace{\zeta_1, \ldots, \zeta_1}_{r_1 \text{ times}}, \underbrace{\zeta_2, \ldots, \zeta_2}_{r_2 \text{ times}}, \ldots, \underbrace{\zeta_m, \ldots, \zeta_m}_{r_m \text{ times}} \},
\end{equation}
with $\sum\limits_{i=1}^{m} r_i = n+p+1$.

The B-spline basis functions are $p$-degree piecewise polynomials on the subdivision $\{\zeta_1,\ldots,\zeta_m\}$. A stable way of generating them is using the Cox-de Boor recursion algorithm \cite{LesPiegl}, which receives as inputs $p$ and $\Xi$. Knot multiplicity is an essential ingredient in spline theory since its allows to control the basis smoothness. Indeed, if a breakpoint $\zeta_j$ has multiplicity $r_j$, then the basis functions have at least $\alpha_j := p - r_j$ continuous derivatives at $\zeta_j$. Hence, the maximum multiplicity allowed for $\zeta_j$ is $r_j = p+1$, in this case $\alpha_j = -1$, and the basis is discontinuous at $\zeta_j$. We restrict ourselves to open knot vectors, in this case $r_1 = r_m = p+1$. Note that this implies $n \geq p+1$ and $\alpha_1 = \alpha_m = -1$. The vector $\boldsymbol{\alpha} := \{\alpha_1,\ldots,\alpha_m\}$ collects the basis regularities. Let's define $\boldsymbol{\alpha} - 1 = \{-1, \alpha_2-1,\ldots,\alpha_{m-1}-1,-1\}$, when $\alpha_j \geq 0$ for $2 \leq j \leq m-1$, and $|\boldsymbol{\alpha}| = \min \{\alpha_2,\ldots,\alpha_{m-1}\}$.

As the name basis suggests, the set $\left\lbrace B^p_i \right\rbrace_{i=1}^n$ defines a linearly independent set of functions with all the good properties wanted for analysis purposes \cite{HBCBOOK}. The space of B-splines spanned by them is denoted by,
\begin{equation}
\mathcal{S}^p_{\boldsymbol{\alpha}} := \text{span} \left\lbrace B^p_i \right\rbrace_{i=1}^n.
\end{equation}
For univariate spline spaces, when $p \geq 1$ and $\alpha_j \geq 0$ for $2 \leq j \leq m-1$, the derivative of a spline is a spline too, indeed the derivative is a surjective operator, that is,
\begin{equation}
\left\lbrace \dfrac{d}{dx} u : u \in \mathcal{S}^p_{\boldsymbol{\alpha}} \right\rbrace \equiv \mathcal{S}^{p-1}_{\boldsymbol{\alpha-1}}. \label{eq:surjectivity_ddx}
\end{equation}
Besides those properties, the refinement of the B-spline spaces by knot insertion and degree elevation are relevant features of spline technology. The knot insertion procedure allows the refinement of the B-spline space by introducing new or repeated knots to the knot vector, guaranteeing the nesting of the unrefined and the refined spaces. In the same way, the degree elevation procedure allows the increase in the polynomial degree of the basis while maintaining the nesting of the unrefined and the refined spaces. It is important to note that such procedures do not commute, and the application of degree elevation followed by knot insertion yields a new refinement procedure, called $k$-refinement. For a better elaboration of such procedures, and examples as well, see \cite{HBC2005},\cite{HBCBOOK} and \cite{LesPiegl}.

\subsection{Bivariate B-splines}
Since we focus on $\R^2$, we restrict our presentation to the bivariate case, but higher dimensional settings are straightforward.

Given $p_1, p_2$, $n_1, n_2$, and the knot vector $\Xi_1$ and $\Xi_2$, we construct a univariate B-spline basis in each direction, that is, $\{ B^{p_d}_{i_d,d} \}_{i_d=1}^{n_d}$ for $d = 1, 2$. The bivariate B-spline basis functions are defined by tensor product of the univariate ones as
\begin{equation}
B^{p_1,p_2}_{i_1,i_2} := B^{p_1}_{i_1,1} \otimes B^{p_2}_{i_2,2} , \quad i_1 = 1,\ldots,n_1;~ i_2 = 1,\ldots,n_2.
\end{equation}

The sets of breakpoints $\boldsymbol{\zeta}_{d} = \{\zeta_{1,d},\ldots,\zeta_{m_d,d}\}$ in each direction $d=1,2$ define the mesh
\begin{equation}
\mathcal{M}_h = \{ Q = (\zeta_{i_1,1},\zeta_{i_1+1,1})\times(\zeta_{i_2,2},\zeta_{i_2+1,2}) : 1 \leq i_1 \leq m_1-1, 1 \leq i_2 \leq m_2-1 \},
\end{equation}
called the parametric mesh, on the parametric domain $\widehat{\Omega} = (0,1)^2$.

Using the notation $\boldsymbol{\alpha}_1 = \{\alpha_{1,1},\ldots,\alpha_{m_1,1}\}$ and $\boldsymbol{\alpha}_2 = \{\alpha_{1,2},\ldots,\alpha_{m_2,2}\}$ for the regularity vectors in each direction, the bivariate B-spline space is defined as
\begin{equation}
\mathcal{S}^{p_1,p_2}_{\boldsymbol{\alpha}_1,\boldsymbol{\alpha}_2} \equiv \mathcal{S}^{p_1,p_2}_{\boldsymbol{\alpha}_1,\boldsymbol{\alpha}_2}(\mathcal{M}_h) := \text{span} \left\lbrace B^{p_1,p_2}_{i_1,i_2} \right\rbrace_{i_1,i_2=1}^{n_1,n_2}.
\end{equation}
The global regularity of the space is defined as $\alpha := \min\{|\boldsymbol{\alpha}_1|,|\boldsymbol{\alpha}_2| \}$.

\subsection{Geometrical mapping and the Physical Mesh}
The subscript $h$ on the parametric mesh notation $\mathcal{M}_h$ stands for a global mesh size, indeed for each $Q \in \mathcal{M}_h$ we define $h_Q := \text{diam}(Q)$ and $h := \max_{Q \in \mathcal{M}_h} h_Q$. In order to guarantee theoretical convergence estimates, the mesh $\mathcal{M}_h$ should satisfy a shape-regularity condition \cite{Bazilevs_Est},
\begin{equation}
\lambda^{-1} \leq \dfrac{h_{Q,\min}}{h_Q} \leq \lambda, \qquad \forall Q \in \mathcal{M}_h,
\end{equation}
for constant $\lambda > 0$, where $h_{Q,\min}$ is the length of the smallest edge of $Q$. If the same $\lambda$ holds for a sequence of nested refined meshes $\{ \mathcal{M}_h \}_{h \leq h_0}$, this sequence is said to be locally quasi-uniform. From now on we will be admitting this assumption.

But the great potential of IGA concept stems from the possibility of working on geometries of varied complexities. This is achieved by the introduction of a geometric mapping $\Fn : \widehat{\Omega} \to \Omega$, from the parametric domain $\widehat{\Omega} = (0,1)^2$ to the general physical domain $\Omega$. We assume that $\Fn$ is a piecewise smooth mapping over $\mathcal{M}_h$, with piecewise smooth inverse. Moreover, $\Fn$ is generally given by B-splines or NURBS basis defined on the coarsest mesh $\mathcal{M}_{h_0}$. The advent of IGA concept started with the observation that $\Fn$ is the object provided by many CAD systems.

Note that implicitly we have a notion of a physical mesh. Indeed, the image of a parametric mesh $\mathcal{M}_h$ induces a mesh on the physical domain $\Omega$, generally denoted by $\mathcal{K}_h$. Also, the images of the elements boundaries by $\Fn$ are denoted by $\mathcal{F}_h$, and boundaries that are contained in $\partial\Omega$ defines the boundary mesh, denoted by $\Gamma_h$.

\section{Divergence-conforming B-spline discretization for the Stokes problem}
\label{sec:divconf}
\subsection{The Stokes problem}
The Stokes system in its strong form is
\begin{subequations}
\begin{align}
- \mathrm{div} (2 \nu \nabla^s \un) + \nabla p &= \fn \text{  in  } \Omega, \label{eq:stokes_momentum}\\
\mathrm{div}~\un &= 0 \text{  in  } \Omega, \label{eq:stokes_incompress}
\end{align}
\end{subequations}
with $\Omega \subset \R^2$ a bounded simply connected Lipschitz open set, $\un$ is the flow velocity, $p$ is the pressure, $\nu > 0$ the kinematic viscosity and $\fn$ denotes a body force acting on the fluid. To be well posed the system must be augmented with appropriate boundary conditions. To simplify the presentation we consider here the case of homogeneous Dirichlet boundary condition, that is, the no-slip case $\un = 0$. Then, as usual in the finite element framework, the strong form is recast in a weak formulation given by:

Find $\un \in \mathbf{H}^1_0(\Omega)$ and $p \in L^2_0(\Omega)$ such that
\begin{equation}
a(\un,\vn) + b(\vn,p) + b(\un,q) = (\fn,\vn)_{\mathbf{L}^2(\Omega)}, \label{eq:stokes_weak}
\end{equation}
for all $\vn \in \mathbf{H}^1_0(\Omega)$ and $q \in L^2_0(\Omega)$ where 
\begin{align}
a(\wn,\vn) &= (2\nu \nabla^s \wn, \nabla^s \vn)_{\mathbf{L}^2(\Omega)}, \\
b(\vn,q) &= -(\mathrm{div}~\vn,q)_{L^2(\Omega)}.
\end{align}

From Brezzi theory \cite{Brezzi1974}, it is known that (\ref{eq:stokes_weak}) is an optimality condition for a saddle-point $(\un,p)$ of a Lagrangian functional, and that a solution $(\un,p) \in \mathbf{H}^1_0(\Omega) \times L^2_0(\Omega)$ exists given that the following conditions hold: the continuity of the bilinear forms $a(\cdot,\cdot)$ and $b(\cdot,\cdot)$, the coercivity of $a(\cdot,\cdot)$, and the inf-sup condition
\begin{equation}
\inf_{q \in L^2_0(\Omega),q \neq 0} \sup_{\vn \in \mathbf{H}^1_0(\Omega)} \dfrac{b(\vn,q)}{||\vn||_{\Hn^1({\Omega})} ||q||_{L^2({\Omega})}} \geq \beta,
\end{equation}
with the constant $\beta > 0$.

\subsection{Divergence-conforming B-spline discretization}
In this section, we review the definitions and results of the divergence-conforming spline discretization for Stokes problem as elaborated for general incompressible flows in \cite{Evans2011}, and first appeared in \cite{BuffaRSV2011}.

Assuming the global regularity $\alpha \geq 1$, the discrete velocity space on the parametric domain $\widehat{\Omega}$ is defined as
\begin{equation}
\widehat{\mathcal{V}}_h := \mathcal{S}^{p_1,p_2-1}_{\boldsymbol{\alpha}_1,\boldsymbol{\alpha}_2-1} \times \mathcal{S}^{p_1-1,p_2}_{\boldsymbol{\alpha}_1-1,\boldsymbol{\alpha}_2} 
\end{equation}
and the discrete pressure space on the parametric domain $\widehat{\Omega}$ as
\begin{equation}
\widehat{\mathcal{Q}}_h := \mathcal{S}^{p_1-1,p_2-1}_{\boldsymbol{\alpha}_1-1,\boldsymbol{\alpha}_2-1}.
\end{equation}
Such pair of spaces can be viewed as smooth generalizations of the Raviart-Thomas elements. Indeed, when $\alpha = 0$ the spaces above coincide with the classical Raviart-Thomas elements, but these elements are not $\mathbf{H}^1$-conforming since they are discontinuous, whereas for $\alpha \geq 1$ the spaces defined above are $\mathbf{H}^1$-conforming, which make them appropriate for discretizing Stokes and Navier-Stokes equations. We are adopting the convention that everything referring to parametric space receives a superscript hat.

The discrete velocity space with no-penetration boundary conditions constraints is defined by
\begin{equation}
\widehat{\mathcal{V}}_{0,h} := \left\{ \widehat{\vn}_h \in \widehat{\mathcal{V}}_h :  \widehat{\vn}_h \cdot \widehat{\nn} = 0 \text{ on } \partial\widehat{\Omega} \right\} \subset \mathbf{H}_0(\widehat{\text{div}};\widehat{\Omega}) ,
\end{equation}
where $\widehat{\nn}$ denotes the outward normal to $\partial\widehat{\Omega}$ and $\widehat{\text{div}}$ the divergence operator in parametric coordinates. With this choice for the velocity space a constrained discrete pressure space
\begin{equation}
\widehat{\mathcal{Q}}_{0,h} := \left\{ \widehat{q}_h \in \widehat{\mathcal{Q}}_h : \int_{\widehat{\Omega}} \widehat{q}_h = 0 \right\} \subset L^2_0(\widehat{\Omega}).
\end{equation}
must be defined.

The rationale of such choices of constrained spaces is that, in order to guarantee a divergence-free velocity field that does not conflict with the {\it inf-sup} stability of the velocity-pressure pair, we must guarantee the surjectivity of the divergence operator at the discrete level. Indeed, together with the surjectivity of the derivative between B-spline spaces (\ref{eq:surjectivity_ddx}), it can be easily seen that
\begin{equation}
\widehat{\mathcal{V}}_{0,h} \xrightarrow{\quad\widehat{\mathrm{div}}\quad} \widehat{\mathcal{Q}}_{0,h}
\end{equation}
forms a cochain complex. Then, if we have the incompressibility condition weakly satisfied, that is,
\begin{equation}
(\widehat{\mathrm{div}}~\widehat{\vn}_h,\widehat{q}_h)_{L^2(\Omega)} = 0 \text{ for all } \widehat{q}_h \in \widehat{\mathcal{Q}}_{0,h},
\end{equation}
we can take $\widehat{q}_h = \widehat{\mathrm{div}}~\widehat{\vn}_h$ above, which implies $||\widehat{\mathrm{div}}~\widehat{\vn}_h||_{L^2(\Omega)} = 0$, and since $\widehat{\mathrm{div}}~\widehat{\vn}_h$ is at least continuous, we have that  $\mathrm{div}~\widehat{\vn}_h = 0$ pointwise.

Note, however, that, if the discrete velocity space is constrained by no-slip conditions on $\partial\widehat{\Omega}$, the discrete pressure space would also be constrained on the corners of $\widehat{\Omega}$, which renders a decrease in the accuracy of the pressure approximation. For a complete discussion see \cite{BuffaFS2011}.

With the growing popularity and successful application of Nitsche's method to impose boundary conditions weakly \cite{Bazilevs2007a}, \cite{Bazilevs2007b}, \cite{Bazilevs2010}, \cite{Harari2010}, the above velocity-pressure pair choice will not be a limitation as the imposition of no-slip conditions could be made weakly by augmenting the variational formulation with additional terms, as we will present on the next section.


Up to now we worked on the parametric domain $\widehat{\Omega}$, but within the IGA framework is easy to work on the physical domain $\Omega$ by means of the geometric mapping $\Fn$. The definition of the discrete velocity and pressure spaces on the physical domain $\Omega$ are made by appropriate transformations induced by $\Fn$. These transformations are a consequence of the pullback operators:
\begin{align}
\iota_\un (\vn) &= \det(D\Fn)(D\Fn)^{-1} (\vn \circ \Fn),  & \vn \in \mathbf{H}_0({\text{div}};{\Omega}) \\
\iota_p   (q)   &= \det(D\Fn)(q \circ \Fn),  & q \in L^2_0({\Omega})
\end{align}
where $D\Fn$ is the Jacobian matrix of the geometrical mapping $\Fn$. The first one, the Piola transform, is a standard choice to build approximation spaces in $\mathbf{H}({\text{div}};{\Omega})$, mainly in the context of mixed finite elements, since it is divergence-preserving and preserves the normal component of the transformed vector field. The second is necessary to preserve the zero mean pressure constraint on the physical domain $\Omega$. 

With the goal of preserving the surjectivity of the divergence operator, but now on physical coordinates, the discrete velocity and pressure spaces on the physical domain are defined by
\begin{align}
\mathcal{V}_{0,h} &:= \left\{  \vn \in \mathbf{H}_0({\text{div}};{\Omega}) :  \widehat{\vn} = \iota_\un (\vn) \in \widehat{\mathcal{V}}_{0,h} \right\}, \\
\mathcal{Q}_{0,h} &:= \left\{  q \in L^2_0({\Omega}) : \widehat{q} = \iota_p (q) \in  \widehat{\mathcal{Q}}_{0,h} \right\}
\end{align}

The last ingredient necessary by the framework of isogeometric differential forms is the existence of suitable projectors. In \cite{BuffaRSV2011}, Buffa et al. introduced suitable $L^2$-stable projection operators $\widehat{\Pi}^0_{\widehat{\mathcal{V}}_{h}} : \mathbf{H}_0({\widehat{\mathrm{div}}};{\widehat{\Omega}}) \to \widehat{\mathcal{V}}_{0,h}$ and $\widehat{\Pi}^0_{\widehat{\mathcal{Q}}_{h}} :  L^2_0(\widehat{\Omega}) \to \widehat{\mathcal{Q}}_{0,h}$ based on dual functionals of B-splines. Then, with the aid of the transformations above, the compositions
\begin{eqnarray}
\Pi^0_{\mathcal{V}_{h}} := \iota_\un^{-1} \circ \widehat{\Pi}^0_{\widehat{\mathcal{V}}_{h}} \circ \iota_\un, \\
\Pi^0_{\mathcal{Q}_{h}} := \iota_p^{-1} \circ \widehat{\Pi}^0_{\widehat{\mathcal{Q}}_{h}} \circ \iota_p,
\end{eqnarray}
render suitable projectors and the main result towards the discrete problem as presented in Evans and Hughes, \cite{Evans2012}:
\begin{proposition}
The diagram
\begin{equation}
\begin{CD}
\mathbf{H}_0({\mathrm{div}};{\Omega}) @>\mathrm{div}>> L^2_0({\Omega}) \\
@V\Pi^0_{\mathcal{V}_{h}}VV  @VV\Pi^0_{\mathcal{Q}_{h}}V \\
\mathcal{V}_{0,h} @>\mathrm{div}>> \mathcal{Q}_{0,h}
\end{CD}
\end{equation}
commutes. Furthermore, there exists a positive constant ${C}_\un$ independent of $h$ such that
\begin{equation}
||{\Pi}^0_{{\mathcal{V}}_{h}} {\vn}||_{\Hn^1({\Omega})} \leq {C}_\un ||{\vn}||_{\Hn^1({\Omega})}, \quad \forall{\vn} \in \mathbf{H}_0({{\mathrm{div}}};{{\Omega}}) \cap \Hn^1({\Omega}).
\end{equation}
\end{proposition}

An important consequence of the proposition above, that follows from the commutativity of the diagram is that, the velocity-pressure pair $( \mathcal{V}_{0,h}, \mathcal{Q}_{0,h} )$ is inf-sup stable. Indeed one can prove:
\begin{proposition}[Evans and Hughes, \cite{Evans2012}]
There exists a positive constant $\beta$ independent of $h$ such that the following holds: for every $q_h \in \mathcal{Q}_{0,h}$, there exists a $\vn_h \in \mathcal{V}_{0,h}$, such that:
\begin{equation}
\mathrm{div} \vn_h = q_h
\end{equation}
and
\begin{equation}
||\vn_h||_{\Hn^1({\Omega})} \leq \beta^{-1} ||q_h||_{L^2({\Omega})}.
\end{equation}
Hence
\begin{equation}
\inf_{q_h \in \mathcal{Q}_{0,h},q_h \neq 0} \sup_{\vn_h \in \mathcal{V}_{0,h}} \dfrac{b(\vn_h,q_h)}{||\vn_h||_{\Hn^1({\Omega})} ||q_h||_{L^2({\Omega})}} \geq \beta
\end{equation}
\end{proposition}

To prove that a velocity-pressure pair is inf-sup stable is generally not an easy task, but surprisingly the proof of the proposition above has seven lines. Moreover, by the surjectivity of the divergence operator (for physical coordinates) the weak satisfaction of the incompressibility condition implies a divergence-free discrete velocity field. 

To maintain compatibility with the notations of \cite{Evans2011},\cite{Evans2012}, in the sections that follow we also define $k' = \min\{p_1-1, p_2-1\}$ to denote the polynomial degree of the velocity-pressure pair. Also, in the numerical examples we must always work with $p = p_1 = p_2$, so in this case, $k'= p-1$ denotes the polynomial degree of the pressure space $\mathcal{Q}_{0,h}$ of the pair.

\section{The discrete variational formulations}
With the discrete divergence-conforming velocity-pressure spaces pair properly defined, we consider the discrete formulation of (\ref{eq:stokes_weak}). Since the discrete velocity space $\mathcal{V}_{0,h}$ only satisfy the no-penetration ($\un \cdot \nn = 0$) constraint, the no-slip condition have to be imposed weakly, that is, by modifying the variational formulation properly. Following \cite{Evans2011}, Nitsche's method is applied. It works as a penalty method by adding variationally consistent terms to the bilinear form $a(\cdot,\cdot)$. Indeed, defining the new bilinear form
\begin{equation}
a_h(\un_h,\vn_h) = a(\un_h,\vn_h) - \eta_h(\un_h,\vn_h),
\end{equation}
where
\begin{equation}
\eta_h(\un_h,\vn_h) = \sum_{F \in \Gamma_h}  \int_F 2 \nu \left( ((\nabla^s \vn_h)~\nn)\cdot\un_h + ((\nabla^s \un_h)~\nn)\cdot\vn_h - \dfrac{C_{pen}}{h_F} \un_h \cdot \vn_h \right) dS, \label{eq:nitsche_lhs}
\end{equation}
the discrete formulation for the no-slip boundary condition is written as:

Find $\un_h \in \mathcal{V}_{0,h}$ and $p_h \in \mathcal{Q}_{0,h}$ such that
\begin{equation}
a_h(\un_h,\vn_h) + b(p_h,\vn_h) + b(q_h,\un_h) = (\fn,\vn_h)_{\mathbf{L}^2(\Omega)}, \label{eq:stokes_discrete}
\end{equation}
for all $\vn_h \in \mathcal{V}_{0,h}$ and $q \in \mathcal{Q}_{0,h}$.

Non-homogeneous tangential Dirichlet boundary conditions are also treated by Nitsche's method. In this case, we add the linear form
\begin{equation}
l_h(\vn_h) = \sum_{F \in \Gamma_h}  \int_F 2 \nu \left( -((\nabla^s \vn_h)~\nn)\cdot\gn + \dfrac{C_{pen}}{h_F} \gn\cdot\vn_h \right) dS \label{eq:nitsche_rhs}
\end{equation}
to the right hand side of (\ref{eq:stokes_discrete}), where $\gn$ is a function defined on $\partial\Omega$ that corresponds to the prescribed tangential component of $\un$ on $\partial\Omega$.

Denoting by $\{ \Phi_1, ..., \Phi_{n_\un} \}$ the basis of $\mathcal{V}_{0,h}$ and by $\{ \phi_1, ..., \phi_{n_p} \}$ the basis of $\mathcal{Q}_{0,h}$, then the solution of (\ref{eq:stokes_discrete}) plus (\ref{eq:nitsche_rhs}) resumes to the solution of the discrete Stokes system
\begin{align}
\begin{bmatrix}
\An_h & B^T \\
B & 0
\end{bmatrix}
\begin{bmatrix}
\un \\
\pn
\end{bmatrix} = 
\begin{bmatrix}
\fn \\
0
\end{bmatrix}
\end{align}
where $\An_h \in \R^{n_{\un} \times n_{\un}}$, $B \in \R^{n_p \times n_{\un}}$ and $\fn \in \R^{n_{\un}}$ are defined by
\begin{align}
&[ \An_h ]_{i,j} = a_h(\Phi_j,\Phi_i) = a(\Phi_j,\Phi_i) - \eta_h(\Phi_j,\Phi_i), \\
&[ B ]_{k,j} = b(\phi_k,\Phi_j), \\
&[\fn]_{i} = (\fn,\Phi_i)_{\mathbf{L}^2(\Omega)} + l_h(\Phi_i),
\end{align}
where $i,j = 1,\ldots,n_{\un}$, $k = 1,\ldots,n_{p}$ and $\un \in \R^{n_{\un}}$ is the coefficient vector of the discrete velocity $\un_h \in \mathcal{V}_{0,h}$ and $\pn \in \R^{n_p}$ is the coefficient vector of the discrete pressure $p_h \in \mathcal{Q}_{0,h}$.

\begin{figure}
\subfigure[Sparsity pattern of Stokes discrete system matrix ($nnz=271662$).]{ 
\includegraphics[scale=0.31]{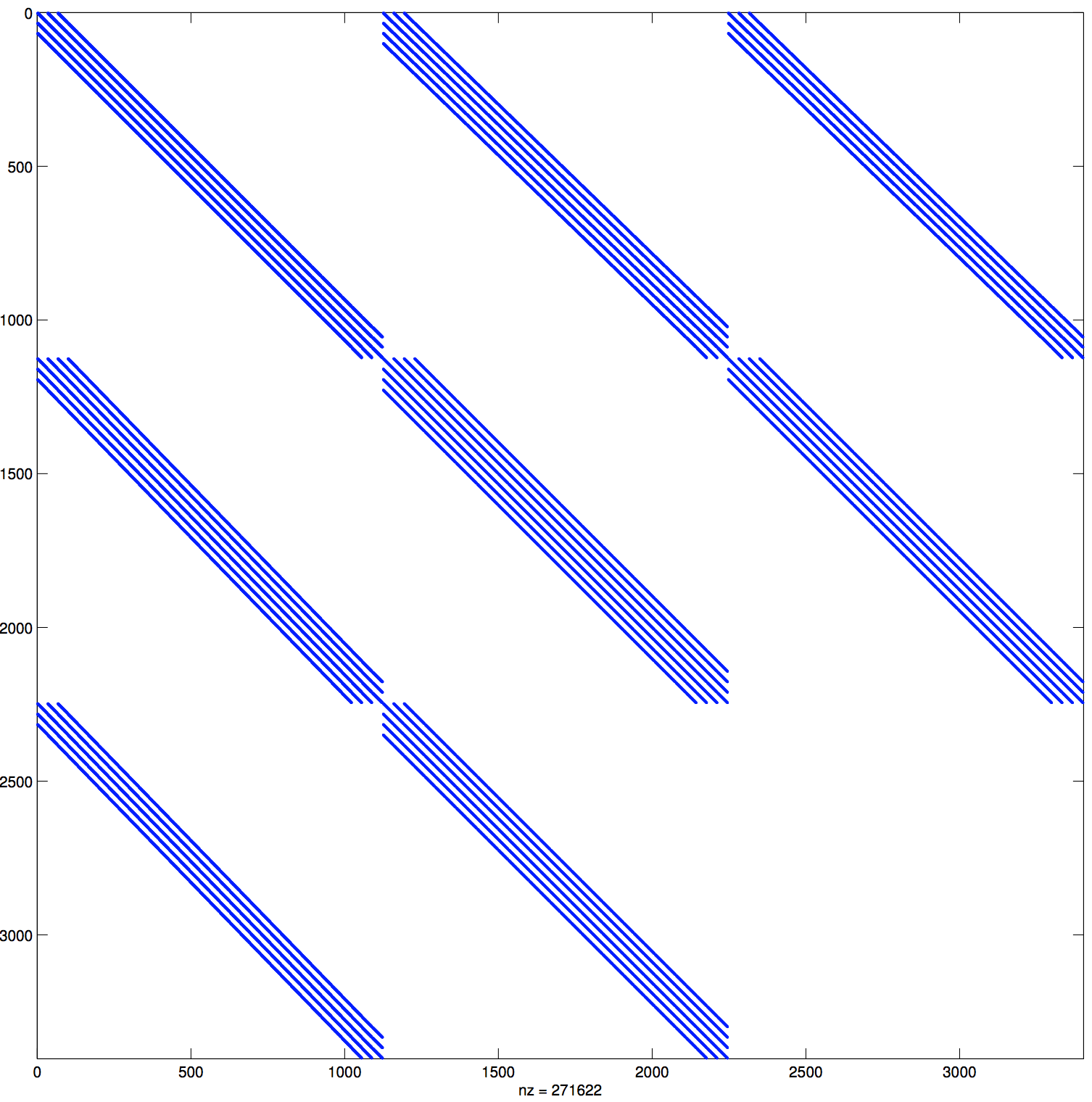} 
\label{fig:spy_mat_deg2_sub32} 
} 
\subfigure[Sparsity pattern of $\An_h$ ($nnz=145634$). The blue dots represent the contributions of $a(\cdot,\cdot)$ and the red dots the contributions of $\eta_h(\cdot,\cdot)$.]{ 
\includegraphics[scale=0.36]{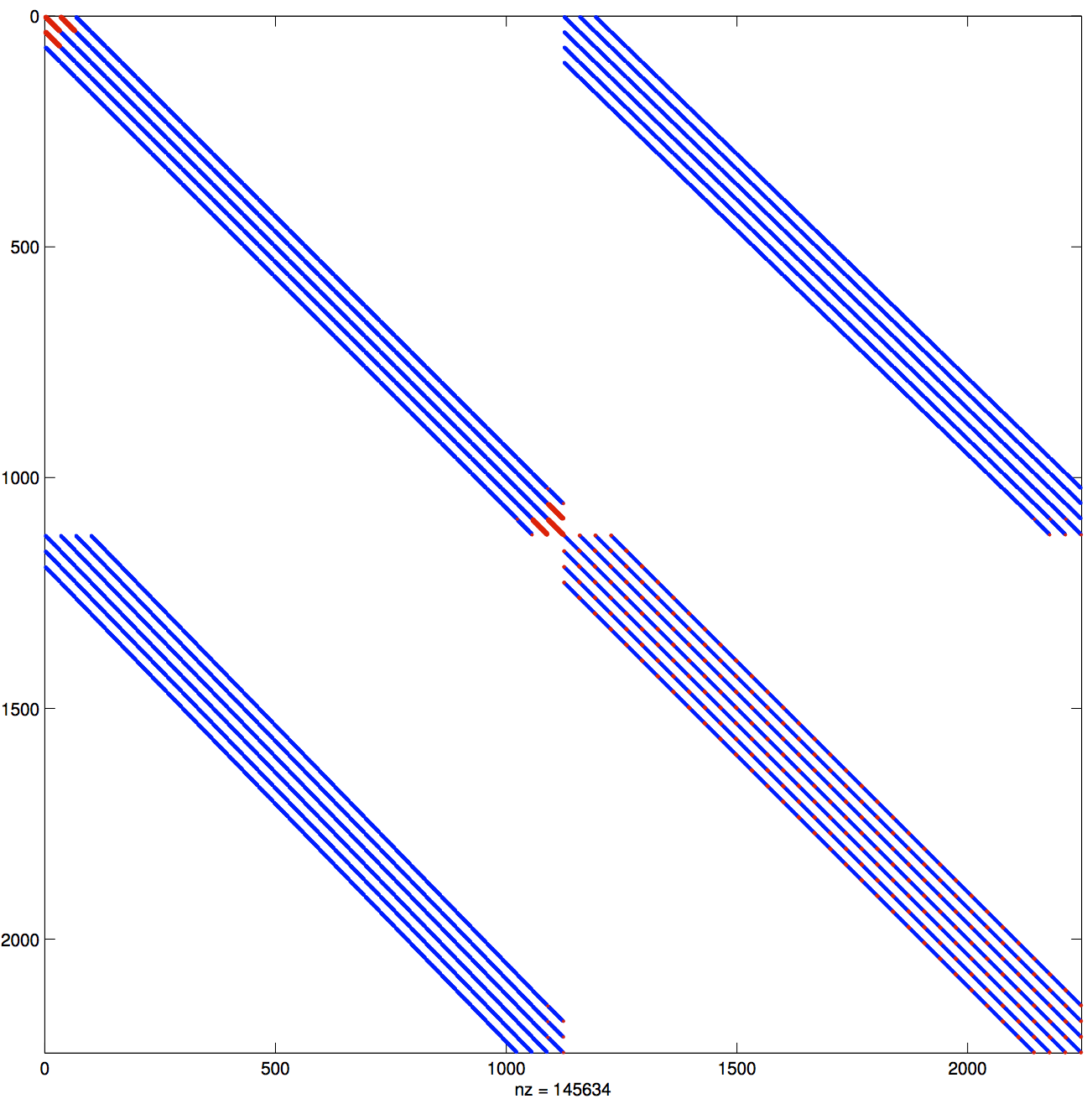} 
\label{fig:spy_A_deg2_sub32} 
}
\caption{Sparsity patterns for $k'=2$, $h=1/32$.}
\label{fig:spy_deg2_sub32}
\end{figure}

\begin{figure}  
\subfigure[Sparsity pattern of Stokes discrete system matrix ($nnz=510990$).]{ 
\includegraphics[scale=0.31]{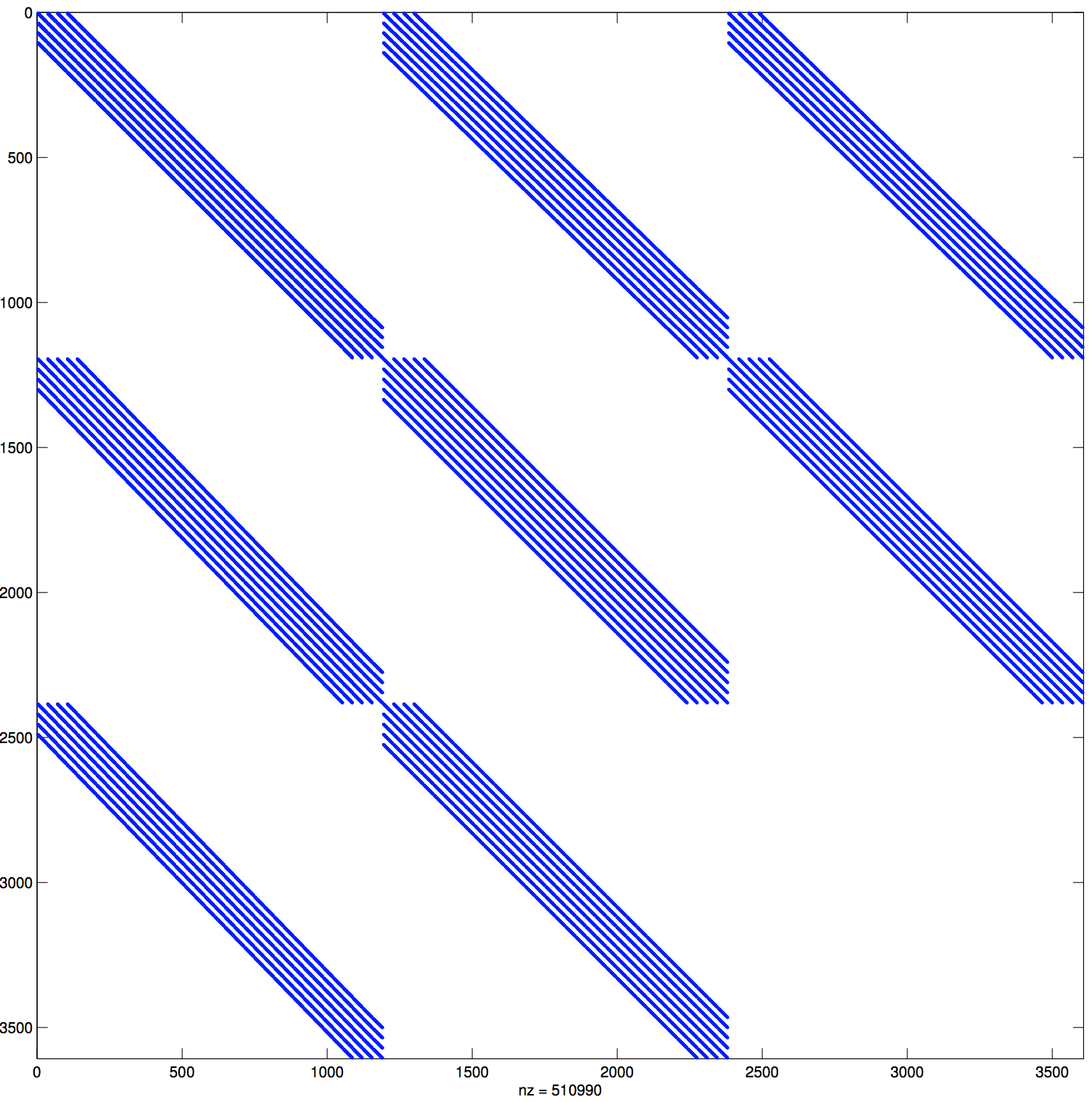} 
\label{fig:spy_mat_deg3_sub32} 
} 
\subfigure[Sparsity pattern of $\An_h$ ($nnz=268606$). The blue dots represent the contributions of $a(\cdot,\cdot)$ and the red dots the contributions of $\eta_h(\cdot,\cdot)$.]{ 
\includegraphics[scale=0.36]{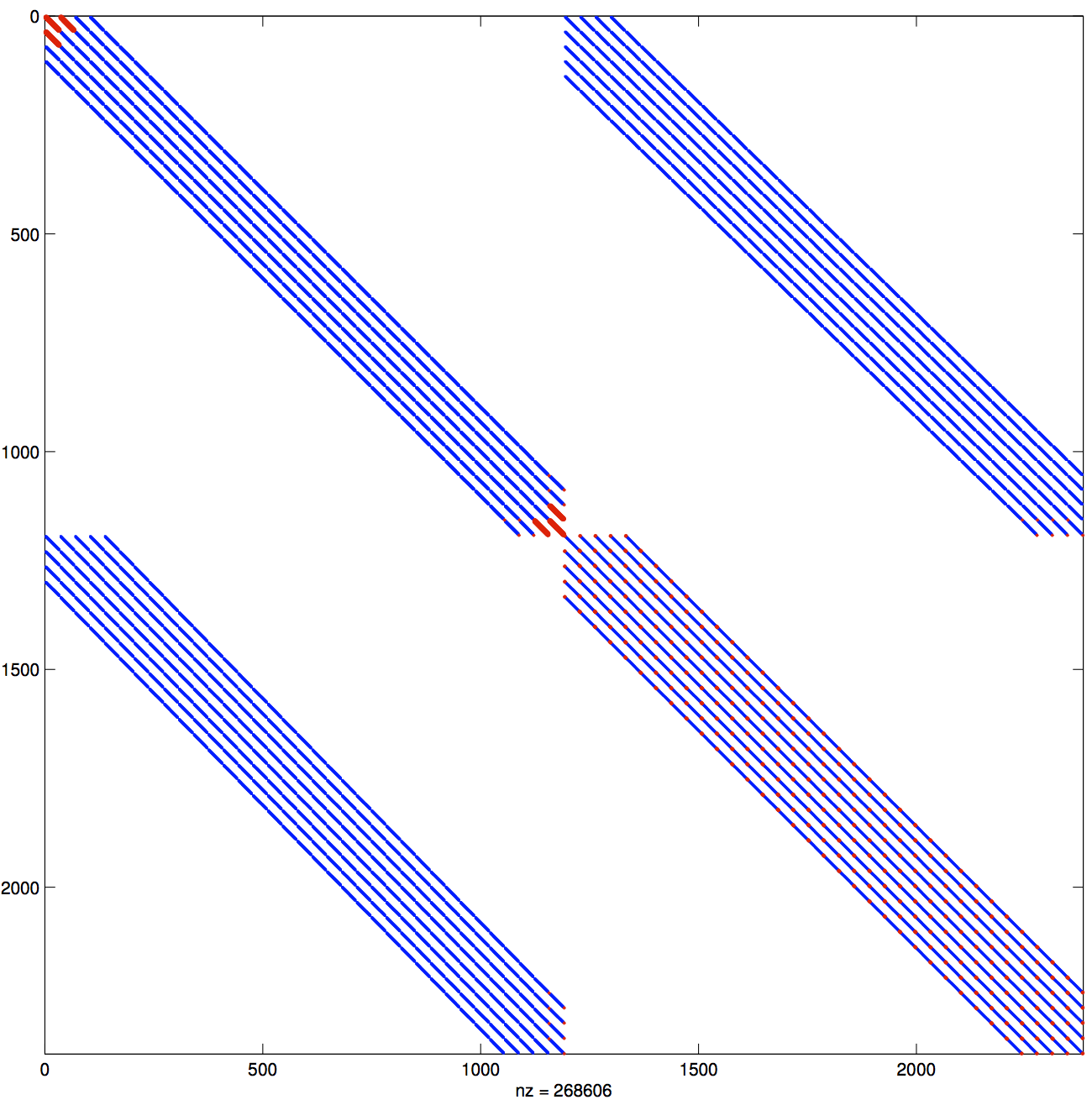} 
\label{fig:spy_A_deg3_sub32} 
}
\caption{Sparsity patterns for $k'=3$, $h=1/32$.}
\label{fig:spy_deg3_sub8}
\end{figure}

Figures \ref{fig:spy_mat_deg2_sub32} and \ref{fig:spy_mat_deg3_sub32} show the sparsity pattern of the coefficient matrix of the discrete Stokes system for the velocity-pressure space pair in question, for the polynomial degrees $k'=2$ and $k'=3$ on $\Omega=(0,1)^2$ for $h=1/32$. Figures \ref{fig:spy_A_deg2_sub32} and \ref{fig:spy_A_deg3_sub32} highlight the (1,1)-block, that is, the matrix $\An_h$, where the blue dots represent the contributions of the bilinear form $a(\cdot,\cdot)$ and the red dots represent the contributions of the bilinear form $\eta_h(\cdot,\cdot)$ for the homogeneous Dirichlet boundary condition problem.

By Lemma 6.3.2 of \cite{Evans2011} we have that, if the $C_{pen}$ is not too small (see Chapter 6,\cite{Evans2011} for an elaboration), the bilinear form $a_h(\cdot,\cdot)$ satisfies
\begin{equation}
\langle \An_h \un, \un \rangle = a_h(\un_h,\un_h) \geq \dfrac{\nu}{C_{Korn}} |\un_h|^2_{\Hn^1(\Omega)} + \sum_{F \in \Gamma_h}  \dfrac{\nu C_{pen}}{2h_F} ||\un_h||^2_{\Ln^2(F)},
\end{equation}
for all $\un_h \in \mathcal{V}_{0,h}$, where $C_{Korn}$ is the constant of Korn's inequality and $|\cdot|_{\Hn^1}$ is the $\Hn^1$-seminorm. In particular, this estimate implies that the symmetric matrix $\An_h$ is positive definite, and that $||\un_h||_{\mathcal{V}} = a_h(\un_h,\un_h)^{1/2} = \langle \An_h \un, \un \rangle^{1/2}$ defines a norm, provided $C_{pen}$ is not too small. 

Combining Lemmas 6.3.1 and 6.3.3 of \cite{Evans2011} we have the stronger discrete {\it inf-sup} stability condition
\begin{equation}
\inf_{p_h \in \mathcal{Q}_{0,h},p_h \neq 0} \sup_{\un_h \in \mathcal{V}_{0,h}} \dfrac{b(p_h,\un_h)}{||\un_h||_{\mathcal{V}} ||p_h||_\mathcal{Q}} \geq \beta_0, \label{eq:stronger_inf_sup}
\end{equation}
where $||p_h||_\mathcal{Q} := \dfrac{1}{2\nu}||p_h||_{L^2(\Omega)}$ and $\beta_0 > 0$ is the stability constant that is independent of $h$, $\nu$, and scales as $O(C_{pen}^{-1/2})$. Tables \ref{tab:inf_sup_penalization_k2} and \ref{tab:inf_sup_penalization_k3} shows the dependency of the stability constant $\beta_0$ on the penalization parameter $C_{pen}$ for $h=1/16$, $\Omega = (0,1)^2$, $k'=2$, and $h=1/16$, $\Omega = (0,1)^2$, $k'=3$, respectively, confirming the dependency numerically.

\begin{table}[h]
\centering
\footnotesize
\begin{tabular}{@{}lrrrrrrrr@{}} 
\toprule 
$C_{pen}$ & $4$ & $8$ & $16$ & $32$ & $64$  & $128$ & $256$ & $1024$\\ 
\midrule 
$\beta_0$ & 0.6015 & 0.5543 & 0.4266 & 0.3106 & 0.2227 & 0.1584 & 0.1122 & 0.0566 \\ 
order & \textemdash &   -0.1178 & -0.3779 & -0.4577 & -0.4801 & -0.4913 & -0.4971 & \textemdash \\ 
\bottomrule 
\end{tabular} 
\caption{Inf-sup stability constant $\beta_0$ for $h=1/16$, $\Omega = (0,1)^2$, $k'=2$ and $\nu=1$.} 
\label{tab:inf_sup_penalization_k2} 
\end{table}

\begin{table}[h]
\centering
\footnotesize
\begin{tabular}{@{}lrrrrrrrr@{}} 
\toprule 
$C_{pen}$ & $4$ & $8$ & $16$ & $32$ & $64$  & $128$ & $256$ & $1024$\\ 
\midrule 
$\beta_0$ & 0.6001 & 0.5853 & 0.4759 & 0.3445 & 0.2456 & 0.1744 & 0.1233 & 0.0616 \\ 
order & \textemdash &  -0.0359 & -0.2985 & -0.4661 & -0.4885 & -0.4940 & -0.5000 & \textemdash \\ 
\bottomrule 
\end{tabular} 
\caption{Inf-sup stability constant $\beta_0$ for $h=1/16$, $\Omega = (0,1)^2$, $k'=3$ and $\nu=1$.} 
\label{tab:inf_sup_penalization_k3} 
\end{table}

Thus, as already remarked by Evans \cite{Evans2011}, the Nitsche's penalization parameter $C_{pen}$ should be taken as small as possible in such a way to guarantee the coercivity of $a_h$, and as we will see in the next section, the inf-sup stability constant $\beta_0$, indeed $\beta_0$ squared, also plays a role in the convergence analysis of the preconditioned MINRES, and it is also desirable to keep it as large as possible (keeping $C_{pen}$ small) for numerical reasons. Finally, in \cite{Evans2011} and \cite{Evans2012} one can find results for stability, existence, uniqueness of the discrete solution, as well as, the mathematical theory of a priori error estimates for the generalized Stokes problem.

\section{Linear Solver}
Here, we discuss the solution of the resulting Stokes system discretized by inf-sup stable mixed elements, that is, the solution of the symmetric indefinite linear system with a block coefficient matrix
\begin{align}
\begin{bmatrix}
\An & B^T \\
B & 0
\end{bmatrix}
\begin{bmatrix}
\un \\
\pn
\end{bmatrix} = 
\begin{bmatrix}
\fn \\
0
\end{bmatrix} \label{eq:stokes_linsys}
\end{align}
where, as usual, $\un \in \R^{n_{\un}}$ is the coefficient vector of the discrete velocity $\un_h$ and $\pn \in \R^{n_p}$ is the coefficient vector of the discrete pressure $p_h$ for the given velocity-pressure pair basis, and $n_{\un}$ is the number of degrees of freedom (number of basis functions) for the velocity and $n_p$ is the number of degrees of freedom for the pressure.

The symmetry property follows from the fact that, generally, $\An$ is a discrete vector Laplacian. The indefinite property follows from Sylvester Law of Inertia \cite{elman2005finite} and the congruence,
\begin{align}
\mathcal{A} = \label{eq:stokes_linsys_lhs}
\begin{bmatrix}
\An & B^T \\
B & 0
\end{bmatrix}
=
\begin{bmatrix}
I & 0 \\
B\An^{-1} & I
\end{bmatrix}
\begin{bmatrix}
\An & 0 \\
0 & -B\An^{-1}B^T
\end{bmatrix}
\begin{bmatrix}
I & \An^{-1}B^T \\
0 & I
\end{bmatrix}.
\end{align}
Indeed, from the congruence above we conclude that $\mathcal{A}$ has $n_{\un}$ positive eigenvalues ($\An$ is positive definite), and it has at least $n_p - 1$ negative eigenvalues, since, for enclosed flows, we have $\mathbf{1} \in Ker(B^T)$, where $\mathbf{1} \in \R^{n_p}$ is the vector with all the components equal to one. The positive semidefinite matrix $S=B\An^{-1}B^T$ is called the {\it pressure Schur complement} and is fundamental in devising good preconditioners.

Despite being symmetric, the indefinite property of $\mathcal{A}$ precludes the application of the Conjugate Gradient method to solve (\ref{eq:stokes_linsys}). Another method from the Krylov subspace methods is better suited for this task: the Minimum Residual Method (MINRES) \cite{Paige1975},\cite{elman2005finite}, which requires only the symmetry of $\mathcal{A}$. For completeness, we review the main characteristics and the convergence results of MINRES and its preconditioned version.

Consider the linear system $\mathcal{A} \xn = \bn$, with $\mathcal{A}$ symmetric and possibly indefinite. Given an initial guess $\xn_0$, the MINRES method generates a sequence of approximate solutions $\xn_k$, $k=1,2,\ldots$ with the property
\begin{equation}
\xn_k \in \xn_0 + \mathcal{K}_k(\mathcal{A},\rn_0),
\end{equation}
where $\rn_0 = \bn - \mathcal{A}\xn_0$ is the initial residual and $\mathcal{K}_k(\mathcal{A},\rn_0) \equiv \text{span} \{\rn_0,\mathcal{A}\rn_0,\ldots,\mathcal{A}^{k-1}\rn_0\}$ is the $k$th dimensional Krylov subspace generated by $\mathcal{A}$ and $\rn_0$. The iterate $\xn_k$ of MINRES is defined satisfying the optimality condition
\begin{align}
||\rn_k||_2 &= \min_{\xn \in \xn_0 + \mathcal{K}_k(\mathcal{A},\rn_0)} ||b - \mathcal{A}\xn||_2 = \min_{p \in \Pi_k} ||p(\mathcal{A})\rn_0||_2 \label{eq:opt_condition} \\
&\leq \min_{p \in \Pi_k} \max_{\lambda \in \sigma(\mathcal{A})} |p(\lambda)| ~ ||r_0||_2. \label{eq:minimax_condition}
\end{align}
where $\Pi_k$ is the set of polynomials of degree at most $k$ with $p(0)=1$, and $\sigma(A)$ is the spectrum of $\mathcal{A}$ (see \cite{Benzi2005survey},\cite{elman2005finite}). The inequality (\ref{eq:minimax_condition}) follows from the symmetry of $\mathcal{A}$. Additionally, in terms of implementation, the symmetry ensures the ultimate efficiency goal of a {\it short-term recurrence} to generate an orthogonal basis for the Krylov subspace $\mathcal{K}_k(\mathcal{A},\rn_0)$, which is achieved by the Lanczos method \cite{ALGA}.

It is clear from (\ref{eq:minimax_condition}) the dependence of MINRES convergence on the spectrum $\sigma(\mathcal{A})$. Since, in this case, $\mathcal{A}$ has positive and negative eigenvalues, clustering such eigenvalues is the goal of good preconditioners for (\ref{eq:stokes_linsys}).

In order to preserve the symmetry of $\mathcal{A}$ a preconditioner $\mathcal{M}$ should be symmetric and positive definite. Therefore, the Cholesky decomposition guarantees the factorization $\mathcal{M} = \mathcal{L}\mathcal{L}^T$. Then, solving the left preconditioned system $\mathcal{M}^{-1}\mathcal{A} \xn = \mathcal{M}^{-1}\bn$ is equivalent to solving the symmetric linear system
\begin{equation}
(\mathcal{L}^{-1}\mathcal{A}\mathcal{L}^{-T}) \yn = \mathcal{L}^{-1}\bn, \qquad \yn = \mathcal{L}^{T}\xn. \label{eq:linsys_precond}
\end{equation}
MINRES can be applied to the preconditioned system above, but now the Euclidean norm of the preconditioned residual $\tilde{\rn}$ is related to the residual of the original system by,
\begin{equation}
||\tilde{\rn}_k||_2 = \langle \mathcal{L}^{-1} \rn_k, \mathcal{L}^{-1} \rn_k \rangle^{1/2} = \langle \mathcal{L}^{-T} \mathcal{L}^{-1} \rn_k, \rn_k \rangle^{1/2} = \langle \mathcal{M}^{-1} \rn_k, \rn_k \rangle^{1/2} = ||\rn_k||_{\mathcal{M}^{-1}},
\end{equation}
and the convergence bound (\ref{eq:minimax_condition}),
\begin{equation}
||\tilde{\rn}_k||_2 \leq \min_{p \in \Pi_k} \max_{\lambda \in \sigma(\mathcal{L}^{-1}\mathcal{A}\mathcal{L}^{-T})} |p(\lambda)| ~ ||\tilde{\rn}_0||_2,
\end{equation}
for the preconditioned system, turns to
\begin{equation}
||\rn_k||_{\mathcal{M}^{-1}} \leq \min_{p \in \Pi_k} \max_{\lambda \in \sigma(\mathcal{M}^{-1}\mathcal{A})} |p(\lambda)| ~ ||\rn_0||_{\mathcal{M}^{-1}}, \label{eq:minimax_precond}
\end{equation}
for the original system. 

We see that the convergence property of the preconditioned MINRES not only depends on $\sigma(\mathcal{M}^{-1}\mathcal{A})$, but is measured in a norm induced by the preconditioner. As observed in \cite{elman2005finite}, the factorization $\mathcal{M} = \mathcal{L}\mathcal{L}^T$ is only needed for theoretical purposes and is never used in practice. Practical implementations only needs the action of $\mathcal{M}^{-1}$, or equivalently, the solution of a linear system with $\mathcal{M}$ as coefficient matrix. Consequently, besides clustering the spectrum of $\mathcal{A}$, a good preconditioner should provide a system of equations that are fast to solve.

\subsection{Block Preconditioning Strategy}
We now discuss a preconditioning strategy for the discrete Stokes system (\ref{eq:stokes_linsys}) \cite{Wathen1993}, \cite{Wathen1994}, \cite{elman2005finite}, \cite{Benzi2005survey}. To start with, consider the block factorization (\ref{eq:stokes_linsys_lhs}). A possible preconditioner in this case is the positive definite block diagonal matrix
\begin{equation}
\mathcal{M} =
\begin{bmatrix}
\An & 0 \\
0 & S
\end{bmatrix}, \label{eq:optimal_precond_M}
\end{equation}
with $S=B\An^{-1}B^T$ ({\it pressure Schur complement}). Indeed, this preconditioner is optimal from the point of view of convergence since it can be proved (Chapter 6, \cite{elman2005finite}) that $\sigma(\mathcal{M}^{-1}\mathcal{A}) = \{ (1-\sqrt{5})/2,1,(1+\sqrt{5})/2 \}$, and then the minimax polynomial convergence estimate (\ref{eq:minimax_precond}) guarantees the convergence of the preconditioned MINRES to the exact solution after at most three iterations. Clearly this preconditioner does not fulfill the requirement of being easily solvable because the solution of a linear system with the pressure Schur complement is not an easy task since it is generally a dense matrix and we do not have $\An^{-1}$ at hand.

A fundamental concept for deriving good preconditioning strategies with the goal of a theoretical scaling property with respect to increasing system size, or equivalently, reducing mesh size $h$, is that of {\it spectral equivalence}. Two matrices $K_1$ and $K_2$ are said to be spectrally equivalent if there are constants $c,C > 0$, both independent of $h$, such that,
\begin{equation}
c \leq \dfrac{\langle K_1 \xn,\xn \rangle}{\langle K_2 \xn,\xn \rangle} \leq C, \qquad \forall \xn \neq \mathbf{0}.
\end{equation}

For general inf-sup stable and conforming mixed discretization, the discrete inf-sup stability condition and the boundedness of the bilinear form $b(\un,p) = -(\mathrm{div}~\un,p)_{L^2(\Omega)}$ (Chapter 5, \cite{elman2005finite}) imply that
\begin{equation}
\beta_0^2 \leq \dfrac{\langle B\An^{-1}B^T \pn, \pn \rangle}{\langle Q \pn, \pn \rangle} \leq C_{b}^2, \quad \forall \pn \in \R^{n_p} \backslash \mathrm{Ker}(B^T),
\end{equation}
where $\beta_0 > 0$ is the inf-sup constant and $C_{b} > 0$ is the boundedness constant, both independent of $h$. The above inequality implies the spectral equivalence of the pressure Schur complement $S=B\An^{-1}B^T$ and the pressure mass matrix $Q$, that is, the matrix whose coefficients are $[Q]_{i,j} = (\phi_j,\phi_i)_{L^2(\Omega)}$, for $i,j = 1,\ldots,n_p$.

Then a more viable preconditioner, in the case of inf-sup stable and conforming discretization, is
\begin{equation}
\mathcal{M} =
\begin{bmatrix}
\An & 0 \\
0 & Q
\end{bmatrix}. \label{eq:ideal_precond_M}
\end{equation}
One can prove in this case (\cite{elman2005finite},\cite{Wathen1994}) that the spectrum $\sigma(\mathcal{M}^{-1}\mathcal{A})$ is included in the union 
\begin{equation}
\left[ \dfrac{1-\sqrt{1+4C_{b}^2}}{2}, \dfrac{1-\sqrt{1+4\beta_0^2}}{2} \right] \bigcup \left\{1 \right\} \bigcup \left[ \dfrac{1+\sqrt{1+4\beta_0^2}}{2}, \dfrac{1+\sqrt{1+4C_{b}^2}}{2} \right]. \label{eq:inclusion_set_ideal}
\end{equation}
In our case, the divergence-conforming discretization, the stronger inf-sup stability condition (\ref{eq:stronger_inf_sup}) and the boundedness of the bilinear form $b$ imply that
\begin{equation}
\beta_0^2 \leq \dfrac{\langle B\An_h^{-1}B^T \pn, \pn \rangle}{\langle Q_\nu \pn, \pn \rangle} \leq C_{b}^2, \quad \forall \pn \in \R^{n_p} \backslash \mathrm{Ker}(B^T),
\end{equation}
where $Q_{\nu} := \frac{1}{2\nu} Q$ is the properly scaled pressure mass matrix since it takes into account the viscosity parameter $\nu > 0$. Then using the preconditioner (\ref{eq:ideal_precond_M}) with $\An = \An_h$ and $Q = Q_{\nu}$ we have the inclusion (\ref{eq:inclusion_set_ideal}) for the spectrum $\sigma(\mathcal{M}^{-1}\mathcal{A})$. Also remember that, in our case $\beta_0^2 = O(C_{pen}^{-1})$.

When the preconditioning system of equations are solved using (\ref{eq:ideal_precond_M}) up to machine precision, like when using a direct solver for each block, it will be called an {\it ideal preconditioning}. In this case, by the eigenvalues bounds (\ref{eq:inclusion_set_ideal}), the inclusion intervals are independent of the mesh-size parameter $h$, and an invariance on the number of iterations for the preconditioned MINRES to converge to a prescribed tolerance, for fixed $\nu$ and $C_{pen}$, is expected. In the next section, we verify numerically that these bounds are indeed sharp, and we have in some sense an optimal clustering of the eigenvalues of the preconditioned system matrix. 

Solving a linear systems with $\An$ and $Q$ maybe not be an easy task. A more general and practical preconditioning strategy is to consider approximations $\An \approx M_{\An}$ and $Q \approx M_Q$, and the preconditioner
\begin{align}
\mathcal{M} =
\begin{bmatrix}
M_{\An} & 0 \\
0 & M_Q
\end{bmatrix}, \label{eq:general_precond_M}
\end{align}
where $M_{\An} \in \R^{n_\un \times n_\un}$ and $M_Q \in \R^{n_p \times n_p}$ are symmetric and positive-definite. The effectiveness of such strategy is given by the following spectral bounds: let $\gamma_{\An}, \Gamma_{\An} > 0$ and $M_{\An}$ be such that
\begin{equation}
\gamma_{\An} \leq \dfrac{\langle \An \un, \un \rangle}{\langle M_{\An} \un, \un \rangle} \leq \Gamma_{\An}, \quad \forall \un \in \R^{n_\un} \backslash \{0\}, \label{eq:spectral_equiv_A}
\end{equation}
and $\gamma_Q, \Gamma_Q > 0$ and $M_Q$ be such that
\begin{equation}
\gamma_{Q} \leq \dfrac{\langle Q \pn, \pn \rangle}{\langle M_{Q} \pn, \pn \rangle} \leq \Gamma_{Q}, \quad \forall \pn \in \R^{n_p} \backslash \{0\}, \label{eq:spectral_equiv_Q}
\end{equation}
and the eigenvalues bounds for the spectrum $\sigma(\mathcal{M}^{-1}\mathcal{A})$ given by the theorem that follows:

\begin{theorem}[Wathen and Silvester \cite{Wathen1994},\cite{elman2005finite}] \label{thm_wathen}
For an inf-sup stable and conforming mixed discretization and a block diagonal preconditioner $\mathcal{M}$ of the form (\ref{eq:general_precond_M}) satisfying (\ref{eq:spectral_equiv_A}) and (\ref{eq:spectral_equiv_Q}), the eigenvalues of the preconditioned Stokes matrix is contained in the union
\begin{equation}
\left[ \dfrac{\gamma_{\An}-\sqrt{\gamma_{\An}^2+4C_{b}^2\Gamma_{\An}\Gamma_{Q}}}{2}, \dfrac{\gamma_{\An}-\sqrt{\gamma_{\An}^2+4\beta_0^2\gamma_{\An}\gamma_{Q}}}{2} \right] \bigcup \left[ \gamma_{\An}, \dfrac{\Gamma_{\An}+\sqrt{\Gamma_{\An}^2+4C_{b}^2\Gamma_{\An}\Gamma_{Q}}}{2} \right].
\label{eq:eigen_bounds}
\end{equation}
\end{theorem}

Obviously several choices are possible for $M_{\An}$ and $M_Q$, but clearly the most robust ones are those where the spectral bounds (\ref{eq:spectral_equiv_A}) and (\ref{eq:spectral_equiv_Q}) are indeed spectral equivalences, that is, do not depend on the mesh-size parameter $h$.

In our numerical tests we consider several possible combinations for $M_{\An}$ and $M_Q$ that we call:
\begin{itemize}
\item Ideal($\An$,$Q$) {\it preconditioning}, where $M_{\An} = \An$ and $M_Q = Q$ and the preconditioner is solved by a direct solver. In this case, the stronger eigenvalues inclusion estimate (\ref{eq:inclusion_set_ideal}) holds;

\item PCG($\An$,$Q$) {\it preconditioning}, where $M_{\An}$ is an approximation for $\An$ by solving the system with the coefficient matrix $\An$ by a preconditioned conjugate gradient method. The same is done with the preconditioner block $M_Q$;

\item Ideal($\An$), Diag$(Q)$ {\it preconditioning}, where $M_{\An} = \An$ is solved by a direct solver and $M_Q = \mathrm{Diag}(Q)$. For classical Lagrangian finite element bases functions it is known \cite{Wathen1993} that $\mathrm{Diag}(Q)$ is spectrally equivalent to $Q$;

\item PCG($\An$), Diag$(Q)$ {\it preconditioning}. This is a combination of the two choices above;

\item and, Diag$(\An,Q)$ {\it preconditioning}, where $M_{\An} = \mathrm{Diag}(\An)$ and $M_Q = \mathrm{Diag}(Q)$. For classical Lagrangian finite element bases functions it is known \cite{Wathen1993} that when $\An$ is a discrete vector Laplacian the spectral bound (\ref{eq:spectral_equiv_A}) holds with $\gamma_{\An} = O(h^2)$ and $\Gamma_{\An} = O(1)$.

\end{itemize}

\section{Numerical Results}
First we present the modifications that had to be done on the MATLAB\textregistered  \textbackslash Octave toolbox GeoPDEs \cite{deFalco2011} to have the divergence-conforming discretization and the Nitsche's method implemented. Then we show the results for three test cases: two manufactured solutions in two different geometries, a square and an 1/8th of an annulus, and the lid-driven cavity flow benchmark.

\subsection{Implementation details}
We discuss some details of GeoPDEs' implementation and some modifications needed to make Nitsche's method work. The main classes are {\bf msh} and {\bf space}. GeoPDE already have an implementation of discretized vector fields mapped by the Piola transformation on the physical domain that is the class {\bf sp\_vector\_2d\_piola\_transform}.

The minor change needed is the introduction of an integral preserving transformation for the discrete pressure space adapted from the class {\bf sp\_bspline\_2d}, renamed to {\bf sp\_scalar\_2d\_integral\_transform}.

The integrals on the boundary are computed also manipulating the classes {\bf msh} and {\bf space}. Both classes have a structure array attribute, called {\bf boundary}, with quadrature information and discrete spaces definitions on the sides of the boundary, respectively. The method {\it msh\_eval\_boundary\_side} of the class {\bf msh2d} already computes the unitary outward normal vector (at the quadrature nodes) on the physical boundary by the expression
\begin{equation}
\nn(\xn) = \dfrac{(D\Fn(\widehat{\xn}))^{-T} \widehat{\nn}(\widehat{\xn})}{||(D\Fn(\widehat{\xn}))^{-T} \widehat{\nn}(\widehat{\xn})||},
\end{equation}
where $\xn = \Fn(\widehat{\xn})$. Since we also need a characteristic length $h_F$ of the (physical) boundary elements we added a new field to the side boundary structure storing the characteristic length of such elements computed as in \cite{Bazilevs2007a}. We augmented the method {\it sp\_eval\_boundary\_side} of the class {\bf sp\_vector\_2d\_piola\_transform} to compute the symmetric gradients of the basis functions evaluated on the quadrature nodes of the boundary. With all the ingredients properly established we defined a function called {\it op\_bilinear\_nitsche\_bnd} to compute the Nitsche's method bilinear operator, and one called {\it op\_linear\_nitsche\_bnd} to compute the Nitsche's method linear functional for the case of non-homogeneous tangential component on the boundary.

\subsection{Square domain}
This example was used for validation and verification in \cite{BuffaFS2011} and \cite{Evans2011}. Here, we will also use it to validate our implementation on GeoPDEs, and to get insight into the preconditioners performance.

\begin{figure}
\centering
\subfigure[Velocity magnitude and some streamlines of the manufactured solution for the square domain.]{ 
\includegraphics[scale=0.33]{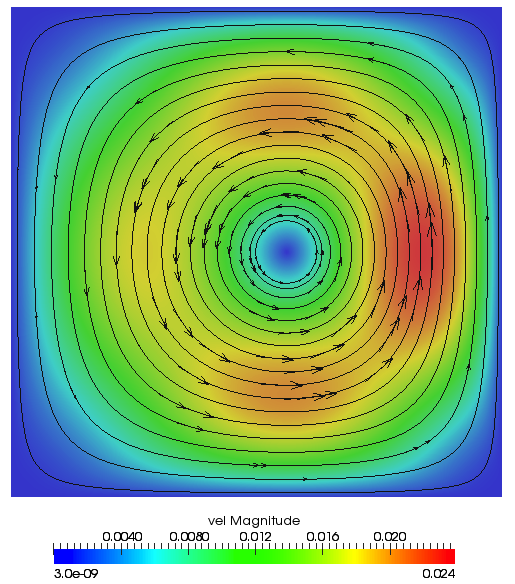}
\label{fig:square_streamlines}
} 
\subfigure[Counter plot of the pressure field of the manufactured solution for the square domain.]{ 
\includegraphics[scale=0.35]{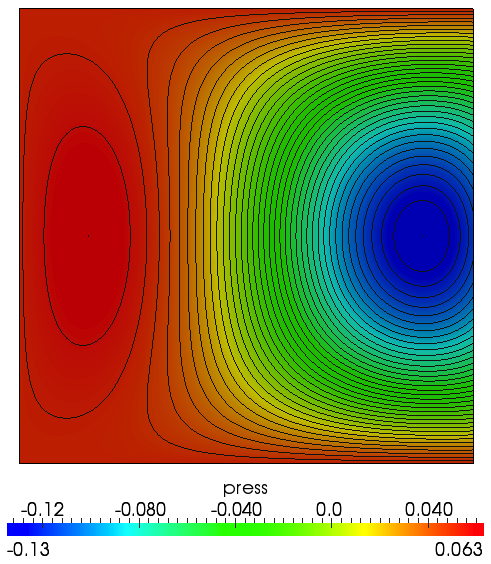}
\label{fig:square_pressure}
}
\caption{Manufactured solution on the square domain.}
\end{figure}

Let $\widehat{\Omega} = \Omega = (0,1)^2$ and the analytical solution
\begin{equation}
\un =
\begin{pmatrix}
2 e^x (-1 + x)^2 x^2 (y^2 - y)(-1 + 2y) \\
e^x (-1 + x) x (-2 + x(3+x)) (-1+y)^2 y^2 
\end{pmatrix}
\end{equation}
and
\begin{align}
p &= -424 + 156e + (y^2 - y) (-456 + e^x (456 + x^2 (228 - 5(y^2 - y)) \\
&+ 2x(-228 + (y^2 - y)) + 2x^3(-36+(y^2 - y)) + x^4(12+(y^2 - y)) ) ).
\end{align}
Figure \ref{fig:square_streamlines} shows the magnitude of the velocity field and some streamlines, whereas Figure \ref{fig:square_pressure} shows the pressure field. The boundary condition is $\un = 0$ on  $\partial\Omega$, the geometric mapping $\Fn$ is the identity mapping and the body force is
\begin{equation}
\fn = - \nabla \cdot (2 \nu \nabla^s \un) + \nabla p. \label{eq:manufactured_force}
\end{equation}

Following \cite{Evans2011}, we choose in all computations the Nitsche's penalization constant as $C_{pen} = 5(k'+1)$. The approximation errors for the velocity and the pressure, and the convergence orders for the polynomials degrees $k'=2$ and $k'=3$ are shown in Tables \ref{tab:square_k2} and \ref{tab:square_k3} respectively. As predicted by the {\it a priori} convergence estimates of \cite{Evans2011}, the order of convergence of the velocity in $\Hn^1$-seminorm is $k'$, and in the $\Ln^2$-norm is $k'+1$. As already observed by Evans \cite{Evans2011}, the computed order of convergence for the pressure in the $L^2$-norm is optimal, that is $k'+1$, whereas the {\it a priori} estimate predicted only a suboptimal order.
\begin{table}[h]
\centering
\footnotesize
\begin{tabular}{@{}lllll@{}} 
\toprule 
$h$ & $1/8$ & $1/16$ & $1/32$ & $1/64$ \\ 
\midrule 
$|\un - \un_h|_{\mathbf{H}^1}$ & 2.26e-03 & 5.58e-04 & 1.39e-04 & 3.46e-05 \\ 
order & \textemdash &   2.02 &   2.01 &   2.00 \\ 
$||\un - \un_h||_{\mathbf{L}^2}$ & 4.10e-05 & 5.18e-06 & 6.55e-07 & 8.26e-08 \\ 
order & \textemdash &   2.98 &   2.98 &   2.99 \\ 
$||p - p_h||_{\mathbf{L}^2}$ & 8.96e-05 & 8.58e-06 & 9.24e-07 & 1.07e-07 \\ 
order & \textemdash &   3.38 &   3.21 &   3.11 \\ 
\bottomrule 
\end{tabular} 
\caption{Errors and convergence orders for $k'=2$ (Square domain).} 
\label{tab:square_k2} 
\end{table} 

\begin{table}[h]
\centering
\footnotesize
\begin{tabular}{@{}lllll@{}} 
\toprule 
$h$ & $1/8$ & $1/16$ & $1/32$ & $1/64$ \\ 
\midrule 
$|\un - \un_h|_{\mathbf{H}^1}$ & 1.24e-04 & 1.62e-05 & 2.09e-06 & 2.65e-07 \\ 
order & \textemdash &   2.93 &   2.96 &   2.98 \\ 
$||\un - \un_h||_{\mathbf{L}^2}$ & 2.33e-06 & 1.58e-07 & 1.03e-08 & 6.53e-10 \\ 
order & \textemdash &   3.88 &   3.95 &   3.97 \\ 
$||p - p_h||_{\mathbf{L}^2}$ & 4.84e-06 & 3.30e-07 & 2.19e-08 & 1.41e-09 \\ 
order & \textemdash &   3.87 &   3.92 &   3.96 \\ 
\bottomrule 
\end{tabular} 
\caption{Errors and convergence orders for $k'=3$ (Square domain).} 
\label{tab:square_k3} 
\end{table}

Tables \ref{tab:niter_square_k_2} and \ref{tab:niter_square_k_3} show the number of iterations and the time in seconds of P-MINRES for all preconditioners discussed in the last section, with a tolerance $tol = 10^{-12}$ for the relative residual. In the cases that PCG was used, the relative residual tolerance is $\sqrt{tol} = 10^{-6}$, and the preconditioner is simply the diagonal of either $\An$ and $Q$. Observe that the almost constant number of iterations of all cases, except Diag($\An,Q$), as the mesh is uniformly refined, indicates that the preconditioning strategies are spectrally equivalent. We decided to postpone a detailed analysis of the preconditioning strategies to the lid-driven cavity flow benchmark.

\begin{table}[h]
\centering 
\footnotesize 
\begin{tabular}{@{}lccccc@{}} 
\toprule 
 & \multicolumn{1}{c}{Ideal($\An,Q$)}   & \multicolumn{1}{c}{PCG($\An,Q$)}   & \multicolumn{1}{c}{Ideal($\An$), Diag($Q$)} &   \multicolumn{1}{c}{PCG($\An$), Diag($Q$)}   & \multicolumn{1}{c}{Diag($\An,Q$)} \\ 
\cmidrule{2-6} 
$h$ & $iter(s)$ & $iter(s)$ & $iter(s)$ & $iter(s)$ & $iter(s)$ \\ 
\midrule 
$1/8$ & 33 (0.79) & 50 (0.28) & 105 (0.33) & 138 (0.40) & 280 (0.07) \\ 
$1/16$ & 35 (0.55) & 45 (0.43) & 125 (1.46) & 152 (0.94) & 526 (0.16) \\ 
$1/32$ & 35 (2.72) & 41 (1.50) & 135 (8.53) & 164 (4.51) & 681 (0.76) \\ 
$1/64$ & 37 (14.17) & 43 (8.67) & 135 (42.60) & 172 (28.43) & 1248 (3.89) \\ 
$1/128$ & 35 (79.74) & 45 (53.05) & 131 (252.68) & 172 (179.45) & 2753 (32.79) \\ 
\bottomrule 
\end{tabular}  
\caption{Number of iterations and times for $k'=2$ (Square domain).} 
\label{tab:niter_square_k_2} 
\end{table}

\begin{table}[h] 
\centering 
\footnotesize 
\begin{tabular}{@{}lccccc@{}} 
\toprule 
 & \multicolumn{1}{c}{Ideal($\An,Q$)}  & \multicolumn{1}{c}{PCG($\An,Q$)}  & \multicolumn{1}{c}{Ideal($\An$), Diag($Q$)}  & \multicolumn{1}{c}{PCG($\An$), Diag($Q$)}  & \multicolumn{1}{c}{Diag($\An,Q$)} \\ 
\cmidrule{2-6} 
$h$ & $iter(s)$ & $iter(s)$ & $iter(s)$ & $iter(s)$ & $iter(s)$ \\ 
\midrule 
$1/8$ & 33 (0.44) & 269 (2.08) & 195 (0.94) & 341 (1.20) & 336 (0.08) \\ 
$1/16$ & 35 (0.90) & 60 (0.71) & 243 (5.14) & 315 (2.50) & 1045 (0.48) \\ 
$1/32$ & 35 (4.49) & 44 (2.46) & 291 (30.92) & 336 (14.53) & 2493 (3.76) \\ 
$1/64$ & 37 (34.27) & 44 (13.02) & 273 (387.19) & 346 (77.33) & 2819 (14.48) \\ 
$1/128$ & 35 (190.75) & 46 (82.38) & 263 (1058.93) & 348 (519.76) & 4431 (83.16) \\ 
\bottomrule 
\end{tabular} 
\caption{Number of iterations and times for $k'=3$ (Square domain).} 
\label{tab:niter_square_k_3} 
\end{table}

\subsection{1/8th annulus domain}
This example, also available at GeoPDEs package \cite{deFalco2011}, is used for validation and verification for the Stokes equations discretizations in \cite{BuffaFS2011}. Here, we will use it with the same goal.

For this example $\Omega \neq \widehat{\Omega}$ is one eighth of an annulus and is parameterized by two types of geometric mappings $\Fn$, the first one is a NURBS parameterization, and the second is a polar parameterization. We used both parameterizations to test if there is any influence of the geometric mapping on the preconditioner and the solver behavior.

\begin{figure}[h] 
\centering 
\includegraphics[scale=0.3]{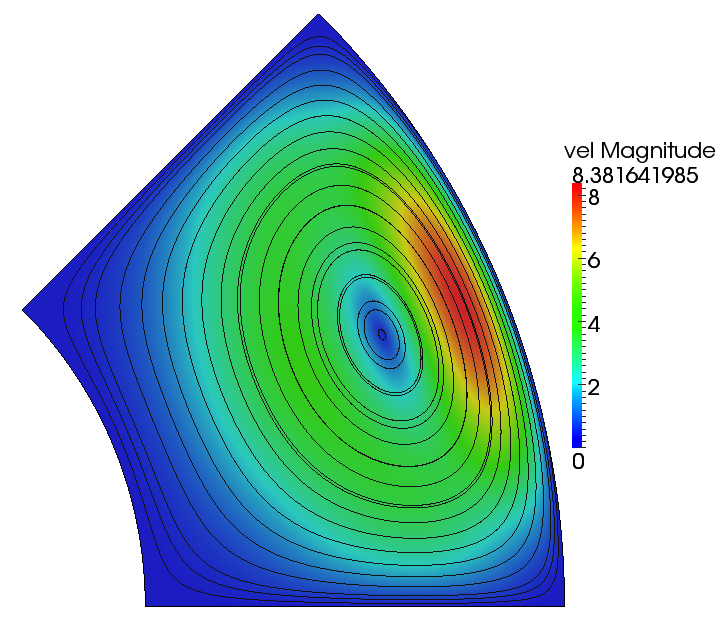} 
\caption{Velocity magnitude and some streamlines of the manufactured solution for the 1/8th annulus domain.}
\label{fig:annulus_streamlines} 
\end{figure}

The boundary condition for this problem is no-slip over $\partial\Omega$, and $\fn$ is given by (\ref{eq:manufactured_force}) for an analytical solution $(\un,p)$ given a priori, that is, a manufactured solution. Figure \ref{fig:annulus_streamlines} shows the domain, the velocity magnitude and some streamlines of the analytical velocity field for this example.

Tables \ref{tab:annulus_errors_k2} and \ref{tab:annulus_errors_k3} show the approximation errors for the velocity and the pressure of the manufactured solution, and the convergence orders for the polynomial degrees $k'=2$ and $k'=3$ respectively, where the geometric mapping $\Fn$ is a NURBS parametrization and for Nitsche's penalization constant $C_{pen} = 5(k'+1)$. Again, as predicted by the {\it a priori} convergence estimates of \cite{Evans2011}, the order of convergence of the velocity in $\Hn^1$-seminorm is $k'$, and in the $\Ln^2$-norm is $k'+1$. As already observed the computed order of convergence for the pressure in the $L^2$-norm is optimal, that is $k'+1$, whereas the {\it a priori} estimate proved in Evans,\cite{Evans2011} is only of suboptimal order. We did not show the approximation errors tables for the case of the polar parameterization since up to the second decimal digit all the errors and convergence orders are the same.
\begin{table}[h] 
\centering
\footnotesize
\begin{tabular}{@{}lllll@{}} 
\toprule 
$h$ & $1/8$ & $1/16$ & $1/32$ & $1/64$ \\ 
\midrule 
$|\un - \un_h|_{\mathbf{H}^1}$ & 3.41e+00 & 8.38e-01 & 2.07e-01 & 5.15e-02 \\ 
order & \textemdash &   2.02 &   2.02 &   2.01 \\ 
$||\un - \un_h||_{\mathbf{L}^2}$ & 6.19e-02 & 7.52e-03 & 9.52e-04 & 1.21e-04 \\ 
order & \textemdash &   3.04 &   2.98 &   2.98 \\ 
$||p - p_h||_{\mathbf{L}^2}$ & 9.00e-02 & 7.71e-03 & 6.62e-04 & 5.75e-05 \\ 
order & \textemdash &   3.55 &   3.54 &   3.53 \\ 
\bottomrule 
\end{tabular} 
\caption{Errors and convergence orders for $k'=2$ (1/8th annulus domain).} 
\label{tab:annulus_errors_k2} 
\end{table} 
\begin{table}[h] 
\centering
\footnotesize
\begin{tabular}{@{}lllll@{}} 
\toprule 
$h$ & $1/8$ & $1/16$ & $1/32$ & $1/64$ \\ 
\midrule 
$|\un - \un_h|_{\mathbf{H}^1}$ & 3.84e-01 & 5.20e-02 & 6.88e-03 & 8.87e-04 \\ 
order & \textemdash &   2.89 &   2.92 &   2.95 \\ 
$||\un - \un_h||_{\mathbf{L}^2}$ & 6.42e-03 & 4.87e-04 & 3.34e-05 & 2.18e-06 \\ 
order & \textemdash &   3.72 &   3.86 &   3.94 \\ 
$||p - p_h||_{\mathbf{L}^2}$ & 3.61e-03 & 1.02e-04 & 4.28e-06 & 2.06e-07 \\ 
order & \textemdash &   5.15 &   4.57 &   4.38 \\ 
\bottomrule 
\end{tabular} 
\caption{Errors and convergence orders for $k'=3$ (1/8th annulus domain).} 
\label{tab:annulus_errors_k3} 
\end{table} 

Tables \ref{tab:niter_annulus_k_2} and \ref{tab:niter_annulus_k_3} show the number of iterations and the time in seconds of P-MINRES for all preconditioners discussed in the last section, with a tolerance $tol = 10^{-12}$ for the relative residual. As in the last example the relative residual tolerance for PCG is $\sqrt{tol} = 10^{-6}$. Similar to the square domain example, the number of iterations of all strategies, except Diag($\An,Q$), are almost constant, indicating the spectral equivalence of these preconditioners. Besides its increasing number of iterations, as mesh is refined, for $k'=2$, Diag($\An,Q$) gave the best time results. Clearly, we do not expect this behavior with more refined meshes. Indeed, for $k'=3$ and for $h \leq 1/32$, Diag($\An,Q$) is slower than PCG($\An,Q$).

\begin{table}[h]
\centering 
\footnotesize 
\begin{tabular}{@{}lccccc@{}} 
\toprule 
 & \multicolumn{1}{c}{Ideal($\An,Q$)}  & \multicolumn{1}{c}{PCG($\An,Q$)}  & \multicolumn{1}{c}{Ideal($\An$), Diag($Q$)}  & \multicolumn{1}{c}{PCG($\An$), Diag($Q$)}  & \multicolumn{1}{c}{Diag($\An,Q$)} \\ 
\cmidrule{2-6} 
$h$ & $iter(s)$ & $iter(s)$ & $iter(s)$ & $iter(s)$ & $iter(s)$ \\ 
\midrule 
$1/8$ & 29 (0.31) & 44 (0.27) & 89 (0.29) & 113 (0.36) & 251 (0.05) \\ 
$1/16$ & 31 (0.47) & 38 (0.38) & 101 (1.22) & 130 (0.82) & 525 (0.17) \\ 
$1/32$ & 31 (2.62) & 38 (1.43) & 105 (6.80) & 140 (3.89) & 728 (0.81) \\ 
$1/64$ & 33 (13.12) & 39 (8.02) & 103 (32.80) & 142 (23.79) & 1216 (4.17) \\ 
$1/128$ & 33 (77.17) & 41 (49.83) & 97 (190.12) & 144 (153.60) & 2814 (33.70) \\ 
\bottomrule 
\end{tabular}  
\caption{Number of iterations and times for $k'=2$ (1/8th annulus domain).} 
\label{tab:niter_annulus_k_2} 
\end{table}

\begin{table}[h]
\centering 
\footnotesize 
\begin{tabular}{@{}lccccc@{}} 
\toprule 
 & \multicolumn{1}{c}{Ideal($\An,Q$)}  & \multicolumn{1}{c}{PCG($\An,Q$)}  & \multicolumn{1}{c}{Ideal($\An$), Diag($Q$)}  & \multicolumn{1}{c}{PCG($\An$), Diag($Q$)}  & \multicolumn{1}{c}{Diag($\An,Q$)} \\ 
\cmidrule{2-6} 
$h$ & $iter(s)$ & $iter(s)$ & $iter(s)$ & $iter(s)$ & $iter(s)$ \\ 
\midrule 
$1/8$ & 31 (0.33) & 69 (0.75) & 113 (0.60) & 203 (0.71) & 341 (0.11) \\ 
$1/16$ & 31 (0.81) & 48 (0.60) & 175 (3.82) & 227 (1.82) & 1044 (0.49) \\ 
$1/32$ & 31 (4.09) & 40 (2.20) & 191 (20.52) & 246 (10.93) & 2328 (3.60) \\ 
$1/64$ & 33 (41.53) & 40 (11.94) & 175 (231.46) & 244 (57.92) & 2882 (14.73) \\ 
$1/128$ & 33 (175.33) & 41 (76.24) & 163 (668.10) & 248 (389.04) & 4371 (82.27) \\ 
\bottomrule 
\end{tabular}
\caption{Number of iterations and times for $k'=3$ (1/8th annulus domain).} 
\label{tab:niter_annulus_k_3} 
\end{table}

\subsection{Lid-driven Cavity Flow}
For the lid-driven cavity flow the parametric domain and the physical domain are the same, that is, the unit square $\widehat{\Omega} = \Omega = (0,1)^2$. The boundary conditions are no-slip at the bottom, and at the sides. At the top the velocity tangential component is constant and in our test equals 1. Finally, the kinematic viscosity is $\nu = 1$.

\begin{figure}
\centering
\subfigure[Velocity magnitude and some streamlines for lid-driven cavity flow.]{ 
\includegraphics[scale=0.4]{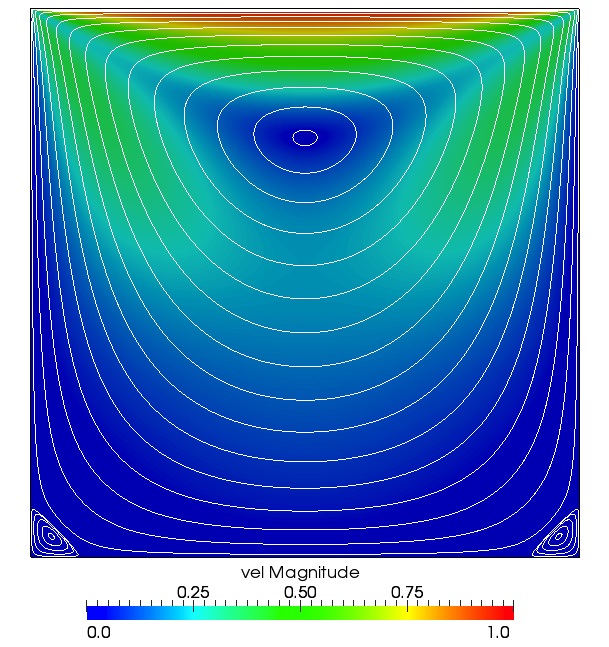}
} 
\subfigure[Zoom in at the left corner showing the primary and the second Moffatt eddies.]{ 
\includegraphics[scale=0.3]{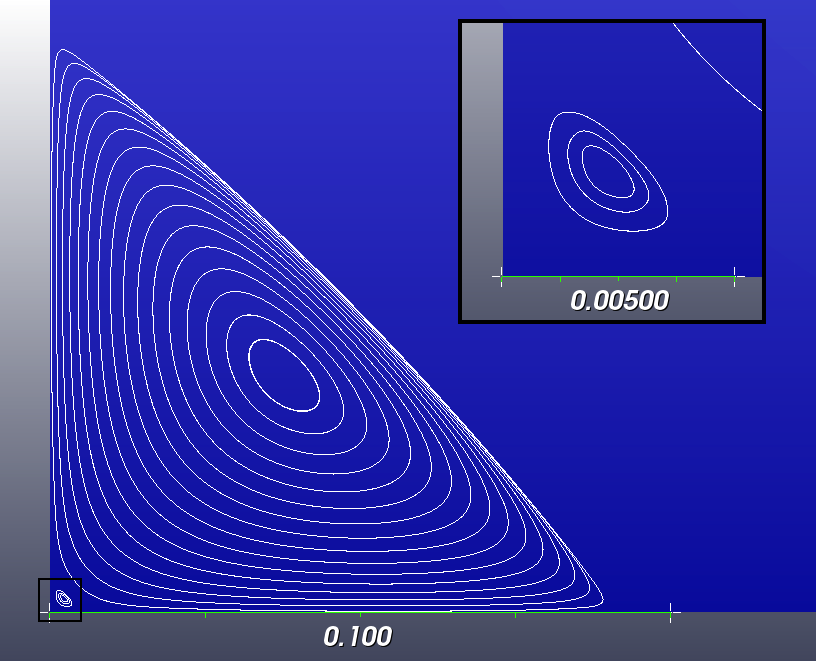}
\label{fig:moffatt_eddy}
}
\caption{Streamlines for the lid-driven cavity flow.}
\label{fig:cavity_streamlines}
\end{figure}

In order the check the correctness of our implementation in figure \ref{fig:cavity_streamlines} we show the streamlines for the lid-driven cavity flow problem. Also, we check an assertion of Evans \cite{Evans2011} that with a uniform refined mesh with mesh-size $h \leq 1/256$, the discretization is able to capture the second Moffatt eddy \cite{Moffatt1964} on the lower corners. Indeed, figure \ref{fig:moffatt_eddy} shows the left lower corner of the domain with the primary and second Moffatt eddies highlighted.


Also, since the cavity flow is an established benchmark problem, we decided to take a closer look on the analysis of the block diagonal preconditioning strategy. To start with, we considered the sizes and numbers of non-zeros of each matrix involved on the linear solver: the coefficient matrix $\mathcal{A}$ of the Stokes system, the viscosity matrix $\An$ and the pressure mass matrix $Q$, for degrees $k'=2$ and $k'=3$, and for five uniform mesh refinements levels. Table \ref{tab:size_nnz_k2} shows the values for the case $k'=2$, and Table \ref{tab:size_nnz_k3} for the case $k'=3$. 

Note that the sizes of the matrices for both cases $k'=2$ and $k'=3$ are almost the same, but the numbers of non-zeros components almost doubles. Obviously this have an impact on the matrix-vector products, but not on vector-vector operations like a dot product or a vector update. A comparison of the cost of matrix-vector operations of a scalar Laplacian for continuous and more regular B-spline bases can be found on Collier et al. \cite{Collier2012}.

\begin{table}[h]
\centering
\footnotesize
\begin{tabular}{@{}ccccccc@{}}
\toprule 
$h$ & $size(\mathcal{A})$ & $nnz(\mathcal{A})$ & $size(\An) = n_{\un}$ & $nnz(\An)$ & $size(Q) = n_p$ & $nnz(Q)$ \\
\midrule
1/8 & 282 & 17654 & 182 & 9170 & 100 & 1936  \\ 
1/16 & 938 & 68774 & 614 & 36482 & 324 & 7056  \\
1/32 & 3402 & 271622 & 2246 & 145634 & 1156 & 26896  \\  
1/64 & 12938 & 1079750 & 8582 & 582050 & 4356 & 104976  \\ 
1/128 & 50442 & 4305734 & 33542 & 2327330 & 16900 & 414736 \\ 
\bottomrule
\end{tabular} 
\caption{Sizes and numbers of non-zeros for the Stokes systems $\mathcal{A}$, $\An$ and $Q$ for $k'=2$.}
\label{tab:size_nnz_k2} 
\end{table}

\begin{table}[h]
\centering
\footnotesize
\begin{tabular}{@{}ccccccc@{}}
\toprule 
$h$ & $size(\mathcal{A})$ & $nnz(\mathcal{A})$ & $size(\An) = n_{\un}$ & $nnz(\An)$ & $size(Q) = n_p$ & $nnz(Q)$ \\
\midrule
1/8 & 343 & 36222 & 222 & 18478 & 121 & 4225  \\ 
1/16 & 1047 & 133294 & 686 & 69342 & 361 & 14641  \\
1/32 & 3607 & 510990 & 2382 & 268606 & 1225 & 54289  \\  
1/64 & 13335 & 2000590 & 8846 & 1057278 & 4489 & 208849  \\ 
1/128 & 51223 & 7916622 & 34062 & 4195198 & 17161 & 819025 \\ 
\bottomrule
\end{tabular} 
\caption{Sizes and numbers of non-zeros for the Stokes systems $\mathcal{A}$, $\An$ and $Q$ for $k'=3$.}
\label{tab:size_nnz_k3} 
\end{table}

Tables \ref{tab:niter_cavity_k_2} and \ref{tab:niter_cavity_k_3} show the numbers of iterations and the times (in seconds) for P-MINRES for the five preconditioning strategies discussed on the last section. The notation of the column PCG($\An,Q$) is as follows: the first number refers to the number of iterations of P-MINRES while the last two is the mean number of iterations for the diagonally preconditioned PCG applied to $\An$ and $Q$, respectively, where the mean is taken with respect to the number of P-MINRES iterations. In all cases, we use a tolerance of $tol = 10^{-12}$ for the relative residual of P-MINRES, and on the cases that we used the diagonally preconditioned PCG the tolerance for the relative residual of PCG is $\sqrt{tol} = 10^{-6}$.

Let's discuss the $k'=2$ case first. Unless the strategy Diag($\An,Q$), all others yield a practically constant number of iterations with respect to mesh refinement. An indication that the preconditioners are spectrally equivalent. Clearly the Ideal($\An,Q$) strategy gave the better results in terms of the number of iterations for the P-MINRES since we solved the preconditioner systems up to machine precision by the backslash command of MATLAB\textregistered, but its time was not the best. Despite the significant increase in the number of iterations, as $h$ is decreased, the Diag($\An,Q$) strategy has the best times for all mesh refinement levels. Obviously we attributed it to the irrelevant cost of solving the preconditioner systems. The second best time is by PCG($\An,Q$).

\begin{table}[h] 
\centering 
\footnotesize 
\begin{tabular}{@{}lccccc@{}} 
\toprule 
 & \multicolumn{1}{c}{Ideal($\An,Q$)}  & \multicolumn{1}{c}{PCG($\An,Q$)}  & \multicolumn{1}{c}{Ideal($\An$), Diag($Q$)}  & \multicolumn{1}{c}{PCG($\An$), Diag($Q$)}  & \multicolumn{1}{c}{Diag($\An,Q$)} \\ 
\cmidrule{2-6} 
$h$ & $iter(s)$ & $iter(s)$ & $iter(s)$ & $iter(s)$ & $iter(s)$ \\ 
\midrule 
$1/8$ & 29 (0.27) & 48/14.29/10 (0.31) & 101 (0.32) & 136 (0.74) & 243 (0.05) \\ 
$1/16$ & 29 (0.46) & 37/25.68/18 (0.36) & 135 (1.59) & 160 (0.99) & 518 (0.16) \\ 
$1/32$ & 29 (2.34) & 37/48.81/34 (1.36) & 135 (9.08) & 166 (4.54) & 696 (0.77) \\ 
$1/64$ & 29 (11.09) & 37/93.65/47.92 (7.51) & 131 (40.97) & 166 (27.16) & 1111 (3.45) \\ 
$1/128$ & 27 (62.15) & 37/183.65/47.19 (43.92) & 123 (237.23) & 164 (171.19) & 2357 (27.55) \\ 
\bottomrule 
\end{tabular} 
\caption{Number of iterations and times for $k'=2$ (lid-driven cavity).} 
\label{tab:niter_cavity_k_2} 
\end{table}

Also, interesting is the worst time, Ideal($\An$), Diag($Q$), followed by PCG($\An$), Diag($Q$). Note also the increase in the number of iterations for these cases as compared to Ideal($\An,Q$) and PCG($\An,Q$) respectively. Our conclusion is that, besides $Q$ and $Diag(Q)$ being spectrally equivalent this is not helping too much P-MINRES. Indeed we numerically computed the spectral bounds $\gamma_Q \approx 0.058$ and $\Gamma_Q \approx 3.33$, where it is clear that the lower bound is considerably small, signaling the deficiency of the diagonal approach. Bellow we will see that it is even worse for the case $k'=3$.

For the Nitsche's penalization constant $C_{pen} = 5(k'+1) = 15$ in this case the squared inf-sup constant is approximately $\beta^2_0 \approx 0.1924$, and the boundedness constant is approximately $C_{b}^2 \approx 1.0134$. Computing the upper bound for the negative eigenvalues given by the inclusion estimate (\ref{eq:eigen_bounds}) of Theorem \ref{thm_wathen} for Ideal($\An$), Diag($Q$) we get $(1 - \sqrt{1 + 4(0.1924)(0.058)})/2 \approx -0.011$. This indicates that the negative part of the spectrum is approaching zero, which is undesirable by the minimax convergence estimate of MINRES, and indeed this happens as values shown in Table \ref{tab:neg_and_pos_eig_k_2}.

\begin{figure}[h]
\centering 
\includegraphics[scale=0.8]{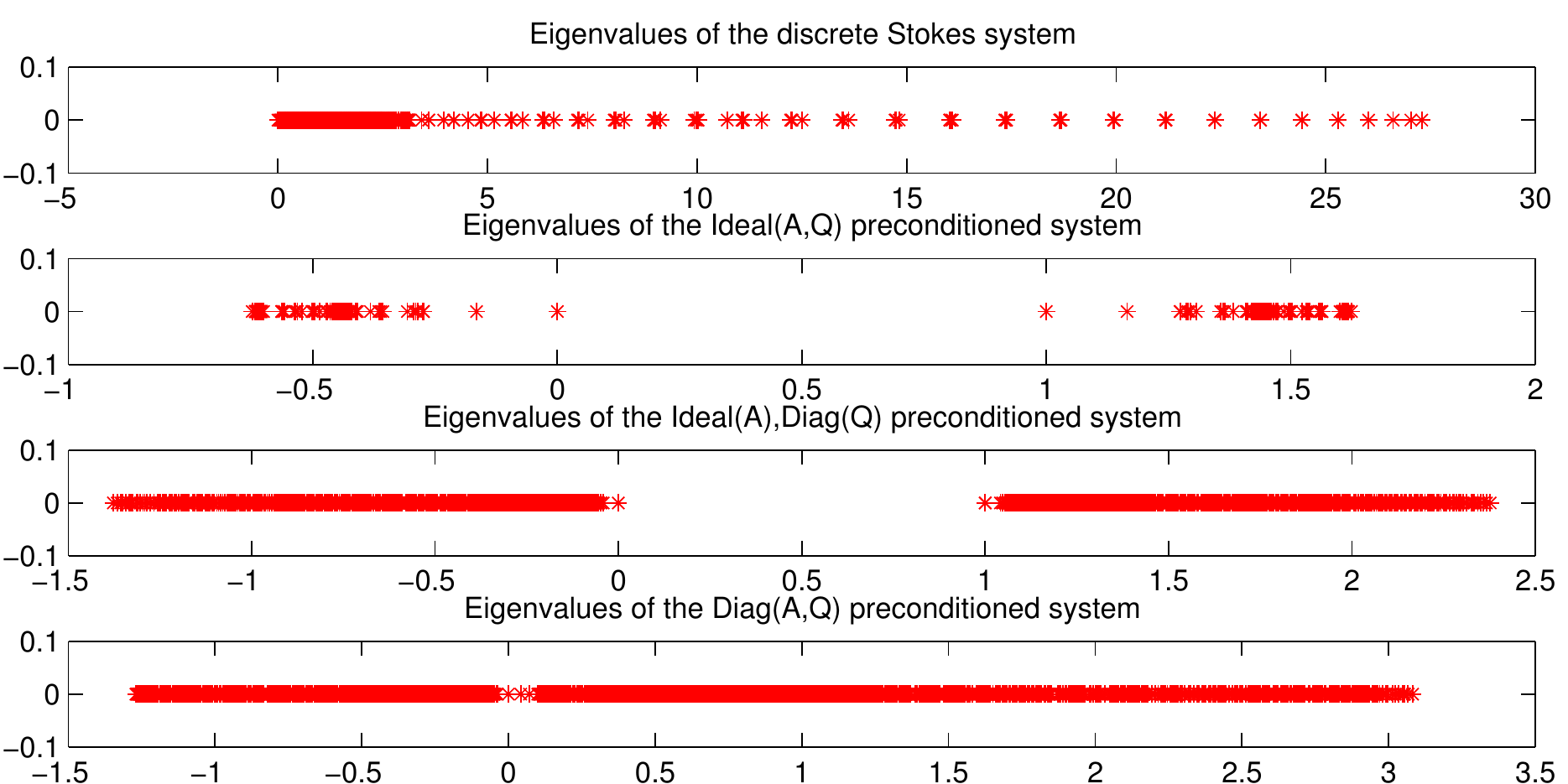} 
\caption{Spectrum of the discrete Stokes system and the preconditioned systems for $k'=2$ and $h=1/32$. (Lid-driven cavity flow)}
\label{fig:spectrum_stokes_cavity_k_2_cpen_15} 
\end{figure}

To corroborate the results of Table \ref{tab:niter_cavity_k_2} we computed numerically $\sigma(\mathcal{A})$ and $\sigma(\mathcal{M}^{-1}\mathcal{A})$ for the preconditioned cases: Ideal($\An,Q$); Ideal($\An$), Diag($Q$) and Diag($\An,Q$). A picture of the spectrum of all cases is shown in figure \ref{fig:spectrum_stokes_cavity_k_2_cpen_15} and some limiting eigenvalues in Table \ref{tab:neg_and_pos_eig_k_2}. We disregard the eigenvalues 0 and 1 because since we imposed the zero mean pressure constraint after the solver, by filtering the solution, the Stokes system matrix is singular and has 0 as an eigenvalue of multiplicity one, that is, the dimension of its kernel. On the other hand, 1 is always an eigenvalue of multiplicity at least $n_{\un} - n_p$ of $\mathcal{M}^{-1}\mathcal{A}$ when $M_{\An} = \An$.

\begin{table}[h]
\centering
\footnotesize
\begin{tabular}{@{}lcccc@{}}
\toprule 
Preconditioner & $\lambda^{-}_{\min}$ & $\lambda^{-}_{\max}$ & $\lambda^{+}_{\min}$ & $\lambda^{+}_{\max}$  \\
\midrule
No precond.  & -4.7522e-4 & -1.3939e-6 & 0.0284 & 27.2898 \\ 
Ideal($\An,Q$) & -0.6240 & -0.1651 & 1.1651 & 1.6240 \\
Ideal($\An$),Diag($Q$) & -1.3773 & -0.0417 & 1.0417 & 2.3773 \\ 
Diag($\An,Q$) & -1.2698 & -0.0375 & 0.0416 & 3.0817 \\  
\bottomrule
\end{tabular} 
\caption{Numerically computed limiting eigenvalues of  $\sigma(\mathcal{M}^{-1}\mathcal{A})$ ($k'=2$).}
\label{tab:neg_and_pos_eig_k_2} 
\end{table}

\begin{figure}[h]
\centering 
\includegraphics[scale=0.6]{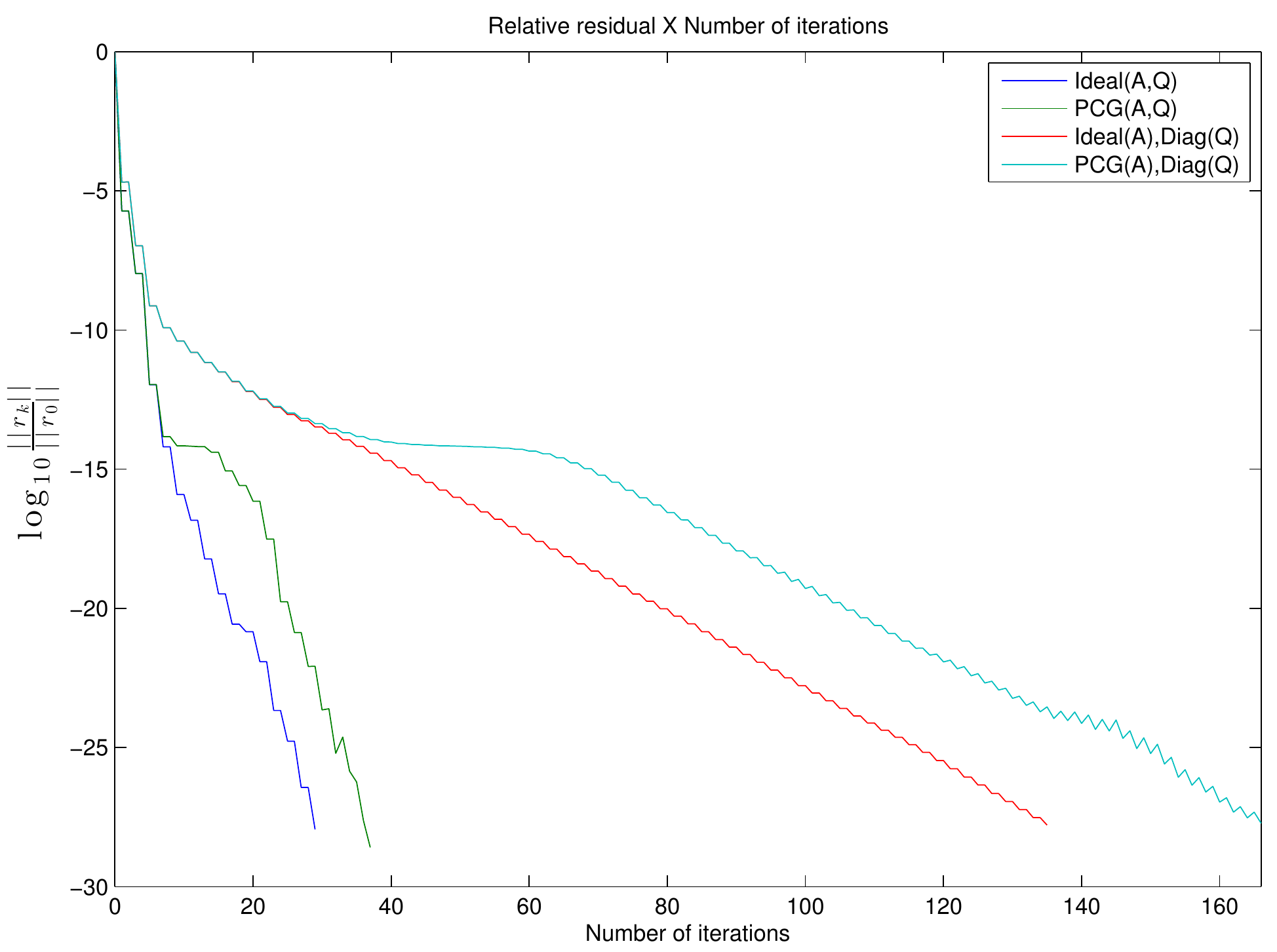} 
\caption{Residuals for P-MINRES $k'=2$ and $h=1/32$. (lid-driven cavity)}
\label{fig:residual_k_2_sub_32} 
\end{figure}

For the cases, Ideal($\An,Q$) and Ideal($\An$), Diag($Q$) one can note that the clustered spectrum attained the symmetry around the value 1/2 as predicted by the inclusion sets estimates (\ref{eq:inclusion_set_ideal}) and (\ref{eq:eigen_bounds}). Also, the bounds for (\ref{eq:inclusion_set_ideal}) are sharp as can be see by computing its values with $\beta^2_0 \approx 0.1924$ and $C_{b}^2 \approx 1.0134$, and comparing with the eigenvalues of Table \ref{tab:neg_and_pos_eig_k_2}. We finish the analysis with figure \ref{fig:residual_k_2_sub_32} showing the relative residual decrease of P-MINRES for the mesh-size $h = 1/32$. As we see the order of convergence for the cases Ideal($\An,Q$) and PCG($\An,Q$) are almost the same as is also for the cases Ideal($\An$), Diag($Q$) and PCG($\An$), Diag($Q$). We conclude that using an iterative solver, as PCG to solve the system for $\An$ does not disturbs the order of convergence of P-MINRES, but clearly using $Diag(Q)$ as an approximation to $Q$ slows down the converge of P-MINRES, as we already concluded by the numerical spectrum shown in figure \ref{fig:spectrum_stokes_cavity_k_2_cpen_15}.

\begin{table}[h] 
\centering 
\footnotesize
\begin{tabular}{@{}lccccc@{}} 
\toprule 
 & \multicolumn{1}{c}{Ideal($\An,Q$)}  & \multicolumn{1}{c}{PCG($\An,Q$)}  & \multicolumn{1}{c}{Ideal($\An$), Diag($Q$)}  & \multicolumn{1}{c}{PCG($\An$), Diag($Q$)}  & \multicolumn{1}{c}{Diag($\An,Q$)} \\ 
\cmidrule{2-6} 
$h$ & $iter(s)$ & $iter(s)$ & $iter(s)$ & $iter(s)$ & $iter(s)$ \\ 
\midrule 
$1/8$ & 31 (0.31) & 92/15.11/9.80 (0.58) & 151 (0.75) & 303 (1.14) & 343 (0.09) \\ 
$1/16$ & 31 (0.80) & 50/27.54/18.98 (0.63) & 277 (5.96) & 325 (2.59) & 1046 (0.49) \\ 
$1/32$ & 31 (4.05) & 38/50.13/35 (2.16) & 315 (33.85) & 364 (15.75) & 2671 (4.04) \\ 
$1/64$ & 29 (57.06) & 38/85.36/67 (10.98) & 297 (451.46) & 366 (82.07) & 2986 (15.37) \\ 
$1/128$ & 29 (167.07) & 38/161.42/116.53 (66.23) & 279 (1121.41) & 362 (521.79) & 4172 (73.64) \\ 
\bottomrule 
\end{tabular} 
\caption{Number of iterations and times for $k'=3$ (lid-driven cavity).} 
\label{tab:niter_cavity_k_3} 
\end{table}

Now we analyze the case $k'=3$. Again in this case we have a numerical indication of the spectral equivalence of the preconditioners, except for the case Diag($\An,Q$). Comparing the columns Ideal($\An,Q$) and PCG($\An,Q$) of Table \ref{tab:niter_cavity_k_3} with the same columns of Table \ref{tab:niter_cavity_k_2} we see that the number of iterations for both discretization degrees are almost the same, but clearly the case $k'=3$ takes more time. This is expected since, for this case, the matrix-vector operations is supposed to take twice the time it takes for the case $k'=2$, because of the increase on the number of non-zeros of the systems matrices for $k'=3$.

In terms of time comparison, we have another picture in this case. Here, the strategy PCG($\An,Q$) has the best times for mesh-sizes $h \leq 1/32$, followed by Diag($\An,Q$). Comparing the number of iterations for Diag($\An,Q$) for both degrees, one can find an increase (for some mesh-sizes more than twice) of the numbers of iterations, that with the additional cost of the matrix-vector operations, led to this result. 

\begin{figure}[h]
\centering 
\includegraphics[scale=0.8]{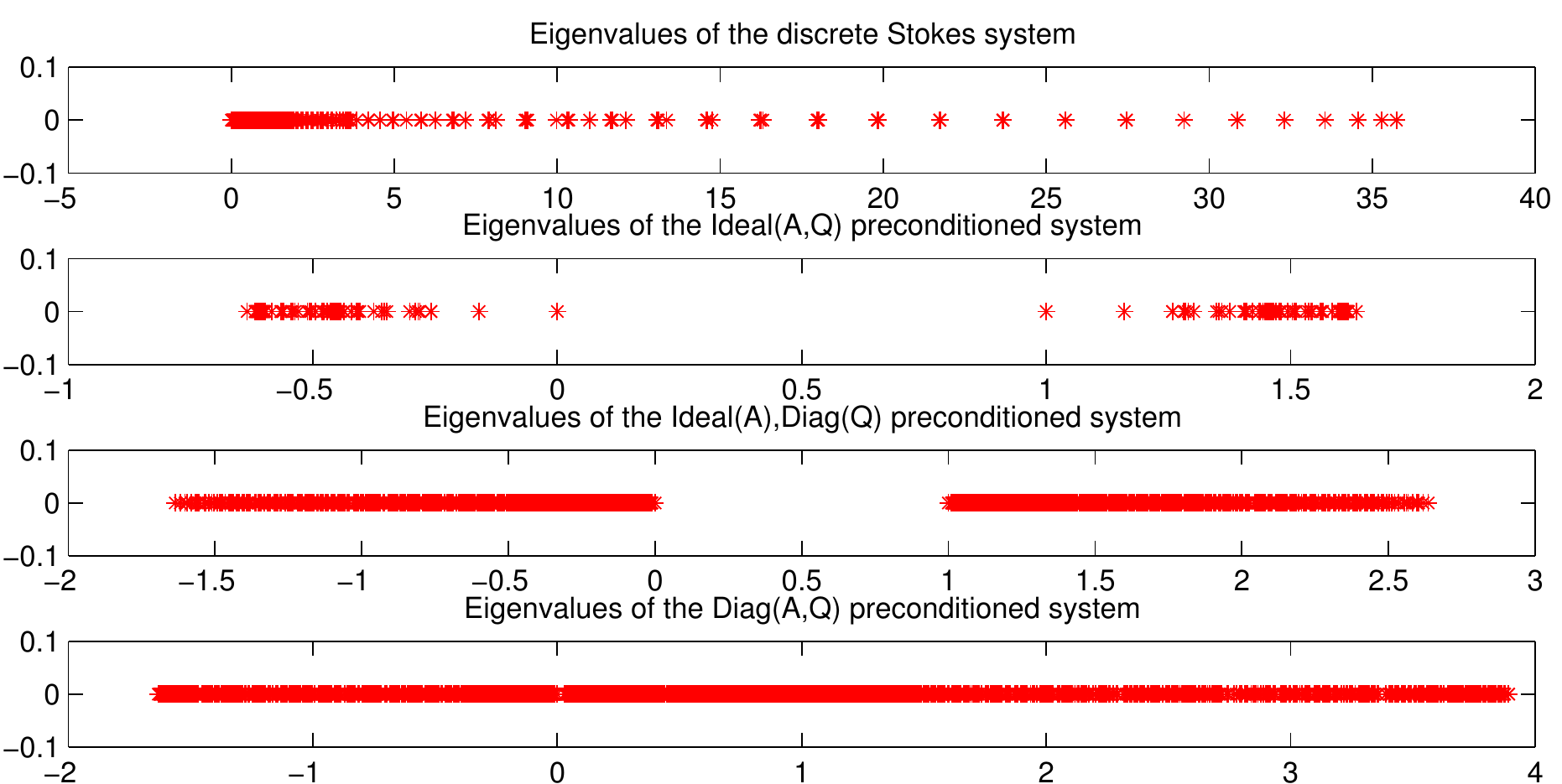} 
\caption{Spectrum of the discrete Stokes system and the preconditioned systems for $k'=3$ and $h=1/32$. (Lid-driven cavity flow)}
\label{fig:spectrum_stokes_cavity_k_3_cpen_20} 
\end{figure}

The worst time continues to be Ideal($\An$), Diag($Q$), followed by PCG($\An$), Diag($Q$), that for most mesh-sizes took half the time of the former, indicating again that besides being spectrally equivalent, Diag($Q$) misses a lot of information about $Q$. Indeed, the numerically computed spectral bounds are $\gamma_Q \approx 0.012$ and $\Gamma_Q \approx 4.46$. It is clear that the lower bound is considerably small, and even smaller than that for $k'=2$ (as we already anticipated), causing the negative part of the spectrum to become closer to zero as is indicated by computing the upper bound $(1 - \sqrt{1 + 4(0.1852)(0.012)})/2 \approx -0.0022$ for the negative eigenvalues given by the inclusion estimate (\ref{eq:eigen_bounds}) of Theorem \ref{thm_wathen} for Ideal($\An$), Diag($Q$), where the squared inf-sup constant is approximately $\beta^2_0 \approx 0.1852$, for the Nitsche's penalization constant $C_{pen} = 5(k'+1) = 20$ in this case (see Figure \ref{fig:spectrum_stokes_cavity_k_3_cpen_20} and Table \ref{tab:neg_and_pos_eig_k_3}). 

\begin{figure}[h]
\centering 
\includegraphics[scale=0.6]{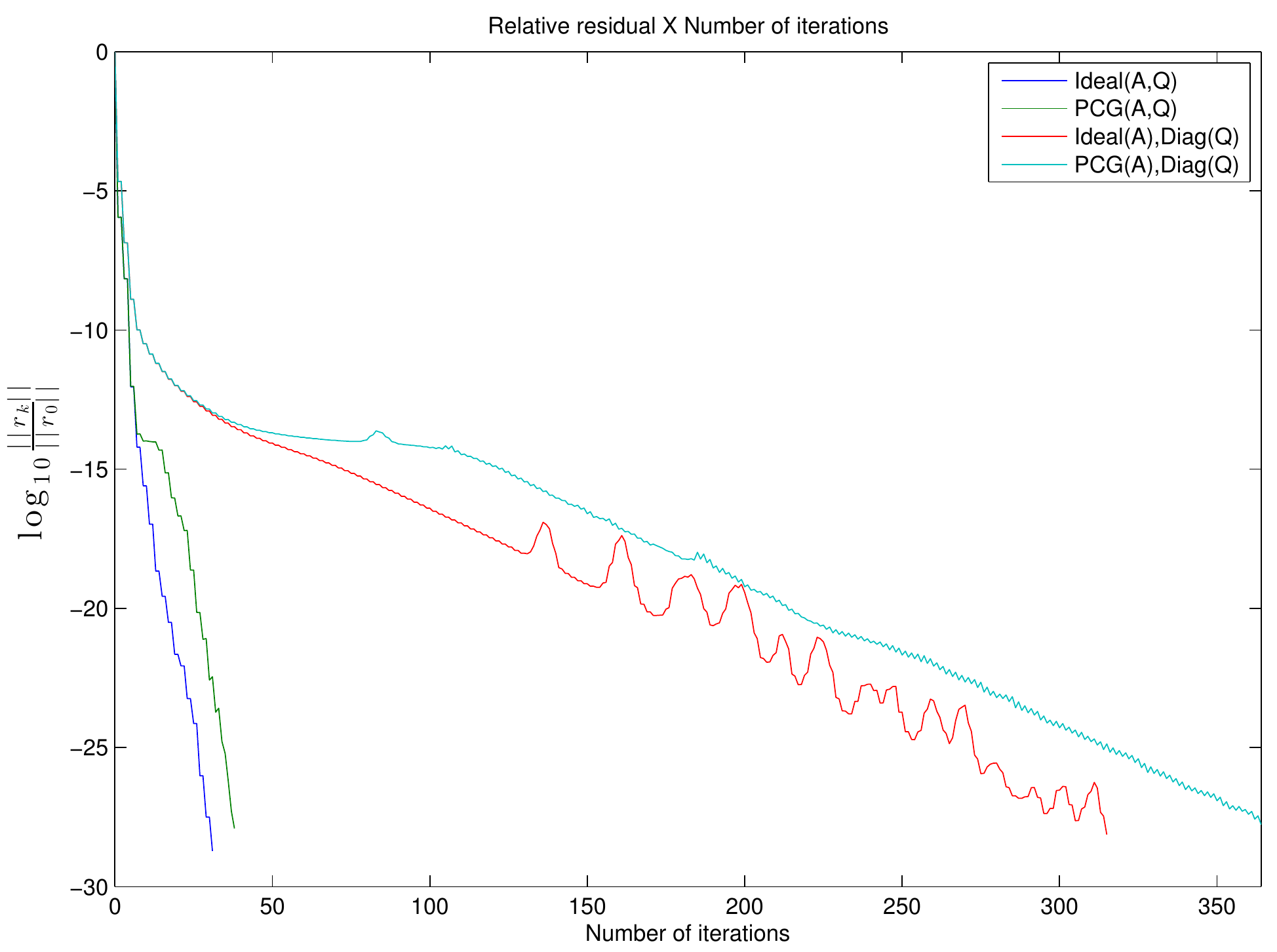} 
\caption{Residuals for P-MINRES $k'=3$ and $h=1/32$. (Lid-driven cavity flow)}
\label{fig:residual_k_3_sub_32} 
\end{figure}

Additionally, it seems by the decrease of the relative residual shown on figure \ref{fig:residual_k_3_sub_32}, that the poor preconditioning offered by $Diag(Q)$, mainly when coupled with Ideal($\An$) causes a non-monotonically convergence of P-MINRES.

\begin{table}[h]
\centering
\footnotesize
\begin{tabular}{@{}lcccc@{}}
\toprule 
Preconditioner & $\lambda^{-}_{\min}$ & $\lambda^{-}_{\max}$ & $\lambda^{+}_{\min}$ & $\lambda^{+}_{\max}$  \\
\midrule
No precond.  & -4.7124e-4 & -3.3468e-7 & 0.0221 & 35.7423 \\ 
Ideal($\An,Q$) & -0.6343 & -0.1597 & 1.1597 & 1.6343 \\
Ideal($\An$),Diag($Q$) & -1.6357 & -0.0107 & 1.0107 & 2.6357 \\ 
Diag($\An,Q$) & -1.6320 & -0.0093 & 0.0299 & 3.8923 \\  
\bottomrule
\end{tabular} 
\caption{Numerically computed limiting eigenvalues of $\sigma(\mathcal{M}^{-1}\mathcal{A})$ ($k'=3$).}
\label{tab:neg_and_pos_eig_k_3} 
\end{table}

For the cases, Ideal($\An,Q$) and Ideal($\An$), Diag($Q$) one can note that the clustered spectrum attained the symmetry around the value 1/2 as predicted by the inclusion sets estimates (\ref{eq:inclusion_set_ideal}) and (\ref{eq:eigen_bounds}). Also, the bounds for (\ref{eq:inclusion_set_ideal}) are sharp as can be see by computing its values with $\beta^2_0 \approx 0.1852$ and $C_{b}^2 \approx 1.0367$, and comparing with the eigenvalues of Table \ref{tab:neg_and_pos_eig_k_3}. We finish the analysis with figure \ref{fig:residual_k_3_sub_32} showing the relative residual decrease of P-MINRES for the mesh-size $h = 1/32$. As we see the order of convergence for the cases Ideal($\An,Q$) and PCG($\An,Q$) are almost the same. The behavior for Ideal($\An$), Diag($Q$) and PCG($\An$), Diag($Q$) are quite different, in spite of appearing that, in the mean, the order of convergence is the same. Up to half of the number of iterations for Ideal($\An$), Diag($Q$) the relative residual decreased monotonically, but after that it started to oscillate up and down. As we concluded above, for $k'=3$ the approximation of $Q$ by $Diag(Q)$ is even worse than that for $k'=2$, and since such oscillations does not appear for Ideal($\An,Q$), we concluded that the discrepancy in the quality of approximations for $\An$ and $Q$ disturbs the convergence of the P-MINRES.

For a more comprehensive comparison, we also tested two global strategies to solve the Stokes system. The first one using also an iterative solver, in this case, the Generalized Minimum Residual Method (GMRES) and the second one using a sparse direct solver, namely the Unsymetric Multifrontal Sparse LU Factorization Package (UMFPACK), called in MATLAB\textregistered by the backslash command. 

We will discuss the iterative solver first. Primarily, a reordering of the unknowns was done using the column approximate minimum degree permutation (COLAMD) algorithm, followed by an ILUT($\tau$) factorization with $\tau = 10^{-6}$ because, in this case, we could not use an ILU(0) factorization since the algorithm breaks down with a zero diagonal element. Then, the factors L and U were used with a left preconditioned GMRES(50), since, in this case, the preconditioner is not symmetric and positive definite. The tolerance used was $10^{-12}$. The results for $k'=2$ and $k'=3$ are shown at Table \ref{tab:gmres_iter}. The total time presented on the last column incorporates the time to setting up the preconditioner, that is, the factorization time, that corresponds approximately to $90 \% $ of the total time.
\begin{table}[h]
\centering
\footnotesize
\begin{tabular}{@{}lcccccr@{}}
\toprule 
$k'$ & $h$ & $nnz(L)$ & $nnz(U)$ & $iter(outer)$ & $iter(inner)$ & $time(s)$ \\
\cmidrule{1-7}
\multirow{2}{*}{2} & 1/64  & 4617086 & 8593117 & 1 & 7 & 21.26 \\ 
& 1/128 & 30400992 & 51347840 & 1 & 9 & 230.02 \\
\cmidrule{1-7}
\multirow{2}{*}{3} & 1/64  & 8380970 & 13280191 & 1 & 9 & 64.84 \\ 
& 1/128 & 47298593 & 77150221 & 1 & 11 & 501.33 \\ 
\bottomrule
\end{tabular} 
\caption{Results for P-GMRES(50) (lid-driven cavity flow).}
\label{tab:gmres_iter} 
\end{table}

The second global strategy, using UMFPACK, gave the time results shown in Table \ref{tab:backslash_solver}. We want to call attention to the case $k'=3$, mesh-size $h=1/256$. In this case, UMFPACK used a standard partial pivoting factorization because the problem is ill-conditioned.
\begin{table}[h]
\centering
\footnotesize
\begin{tabular}{@{}lcr@{}}
\toprule 
$k'$ & $h$ & $time(s)$ \\
\cmidrule{1-3}
\multirow{3}{*}{2} & 1/64  &  1.70 \\ 
                   & 1/128 &  12.94 \\
                   & 1/256 &  435.00 \\
\cmidrule{1-3}
\multirow{3}{*}{3} & 1/64  &  5.77 \\ 
                   & 1/128 &  36.36 \\
                   & 1/256 &  2289.70 \\ 
\bottomrule
\end{tabular} 
\caption{Results for UMFPACK (lid-driven cavity flow).}
\label{tab:backslash_solver} 
\end{table}

For real larger problems, some computational techniques should be used in advance. One such technique is matrix reordering that, in addition to possibly improving data locality, also helps to improve the quality of the preconditioner for pure algebraic strategies, like incomplete factorizations. For a brief review see \cite{Camata2011}, Section 2. Also, Collier et al. \cite{Collier2012} showed that incomplete factorization with zero fill-in performs well as a preconditioner for the conjugate gradient method for isogeometric discretization of a Laplace problem. Moreover, the incomplete factorization preconditioning presented $p$-scalability, that is, scalability under polynomial refinement for $C^{p-1}$ bases, but no spectral equivalence with mesh refinement.

\begin{figure}
\subfigure[Sparsity patterns for RCM of $\An$ ($k'=2$ and $h=1/32$).]{ 
\includegraphics[scale=0.35]{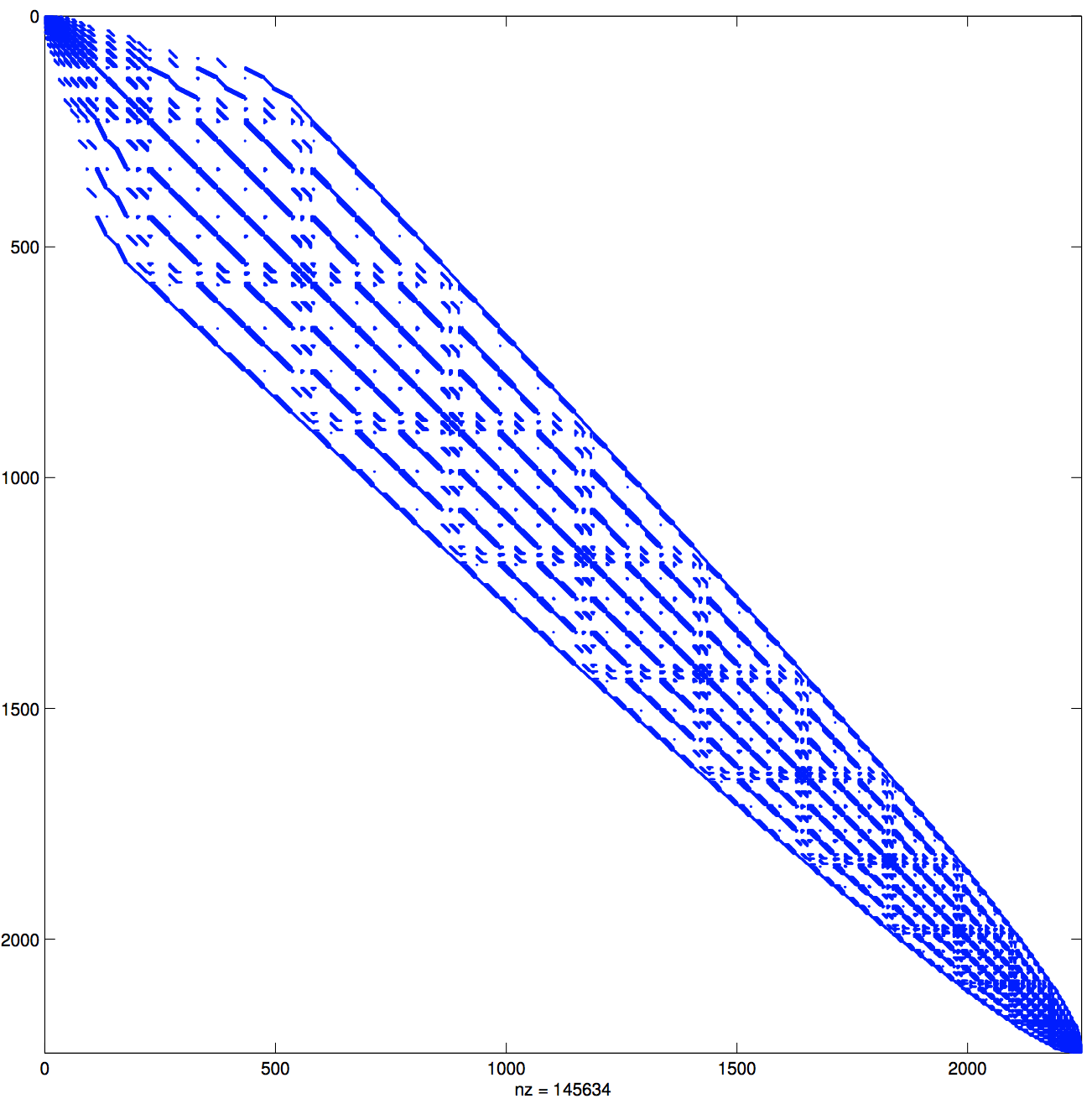} 
} 
\subfigure[Sparsity patterns for RCM of $\An$ ($k'=3$ and $h=1/32$).]{ 
\includegraphics[scale=0.35]{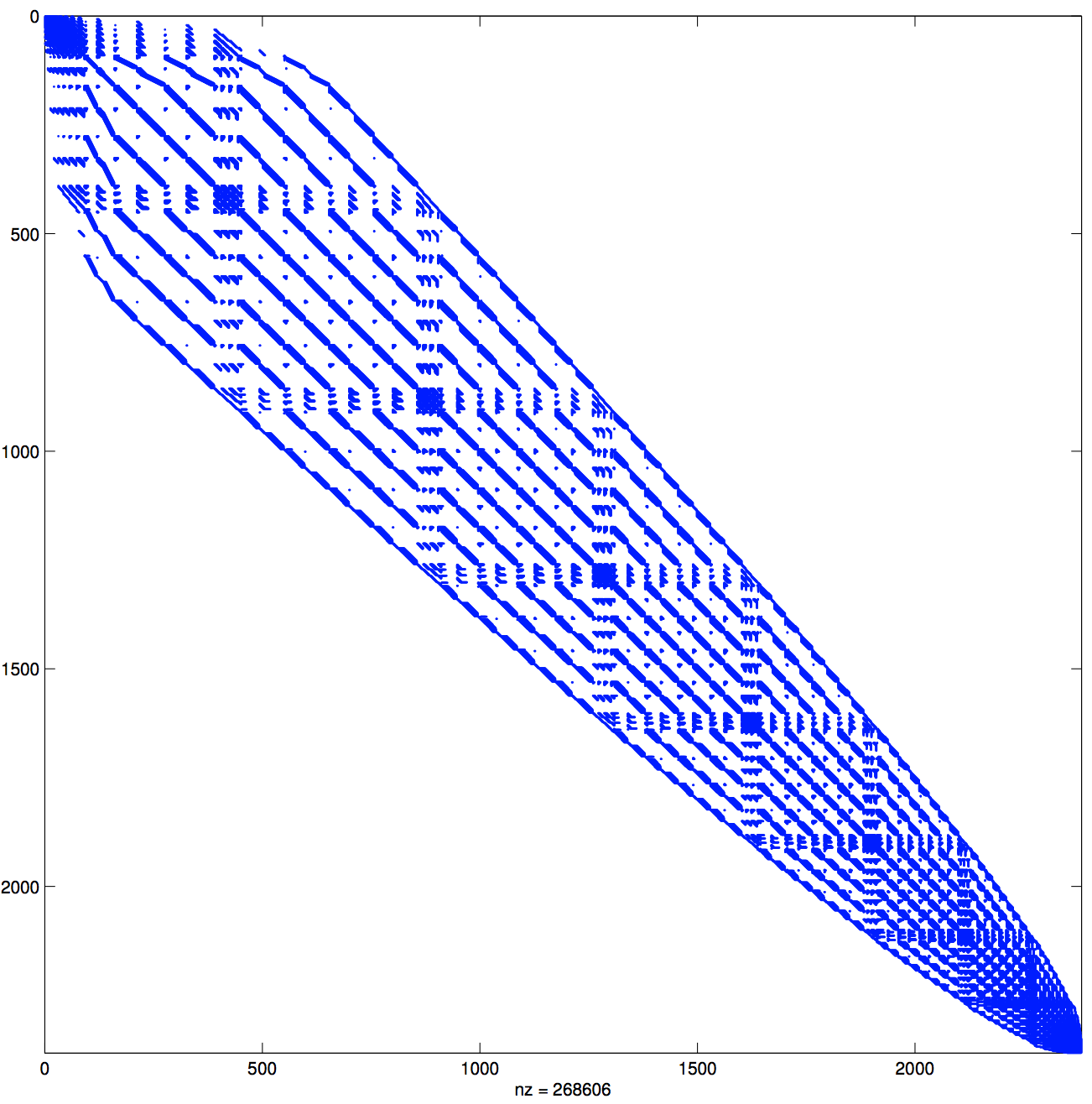} 
}
\caption{Sparsity patterns for RCM of $\An$.}
\label{fig:spy_rcm_sub32}
\end{figure}

In this context, the last numerical experiment is as follows: first we promoted a separate reordering of the matrices $\An$ and $Q$ using the reverse Cuthill-McKee (RCM) algorithm, shown in Figure \ref{fig:spy_rcm_sub32} (compare with figures \ref{fig:spy_A_deg2_sub32} and \ref{fig:spy_A_deg3_sub32}), then we did an incomplete Cholesky factorization with zero fill-in, IC(0), of both $\An$ and $Q$. Finally, the factors were used as preconditioners for the PCG method to solve the preconditioner systems with $\An$ and $Q$ for P-MINRES. We also experimented the effect of a relaxation on the relative residual tolerance of PCG; hence we tested with $tol = 10^{-6}$ as in all tests above, and with $tol = 10^{-3}$. The results are shown in Tables \ref{tab:ic_pgc_tols_k2} and \ref{tab:ic_pgc_tols_k3}.

\begin{table}[h] 
\centering 
\footnotesize
\begin{tabular}{@{}lcccccccc@{}} 
\toprule[1.5pt] 
 & \multicolumn{8}{c}{IC(0)-PCG($\An,Q$)} \\ 
\cmidrule[1pt]{2-9}
 & & & & \multicolumn{2}{c}{$tol=10^{-3}$} & & \multicolumn{2}{c}{$tol=10^{-6}$} \\
\cmidrule{5-6}  \cmidrule{8-9}
$h$ & $R+F(s)$ & $nnz(IC(\An))$ & $nnz(IC(Q))$ & $iter$ & $time(s)$ & & $iter$ & $time(s)$ \\ 
\midrule 
$1/8$ & 9.64e-4 & 4676 & 1018 & 39/2.10/1 & 0.10 & & 34/4.18/1  & 0.17  \\ 
$1/16$ & 0.0033 & 18548 & 3690 & 40/3.35/1 & 0.16 & & 36/6.89/1 & 0.18  \\ 
$1/32$ & 0.0125 & 73940 & 14026 & 42/6.14/1 & 0.51 & & 37/12.65/1 & 0.71  \\ 
$1/64$ & 0.0471 & 295316 & 54666 & 43/11.05/1 & 2.60 & & 36/24.06/1 & 4.16  \\ 
$1/128$ & 0.1824 & 1180436 & 215818 & 44/19.48/1 & 18.11 & & 37/46.70/1 & 33.46 \\
$1/256$ & 0.7136 & 4720148 & 857610 & 42/36.40/1 & 116.10 & & 36/90.92/1 & 237.45 \\
\bottomrule[1.5pt]
\end{tabular}
\caption{Number of iterations of P-MINRES, mean number of iterations for the PCG($\An$) and PCG($Q$), and total time for $k'=2$ (lid-driven cavity).} 
\label{tab:ic_pgc_tols_k2} 
\end{table}

\begin{table}[h] 
\centering 
\footnotesize
\begin{tabular}{@{}lcccccccc@{}} 
\toprule[1.5pt] 
 & \multicolumn{8}{c}{IC(0)-PCG($\An,Q$)} \\ 
\cmidrule[1pt]{2-9}
 & & & & \multicolumn{2}{c}{$tol=10^{-3}$} & & \multicolumn{2}{c}{$tol=10^{-6}$} \\
\cmidrule{5-6}  \cmidrule{8-9}
$h$ & $R+F(s)$ & $nnz(IC(\An))$ & $nnz(IC(Q))$ & $iter$ & $time(s)$ & & $iter$ & $time(s)$ \\ 
\midrule 
$1/8$ & 0.0024 & 9350 & 2173 & 42/1.95/1 & 0.18 & & 34/3.56/1  & 0.18  \\ 
$1/16$ & 0.0074 & 35014 & 7501 & 40/2.55/1 & 0.18 & & 37/5.11/1 & 0.21  \\ 
$1/32$ & 0.0283 & 135494 & 27757 & 40/4.48/1 & 0.61 & & 37/9.14/1 & 0.85  \\ 
$1/64$ & 0.1053 & 533062 & 106669 & 43/7.63/1 & 2.96 & & 37/17.41/1 & 4.68  \\ 
$1/128$ & 0.4038 & 2114630 & 418093 & 44/13.80/1 & 21.06 & & 38/33.66/1 & 37.81 \\
$1/256$ & 1.6763 & 8423494 & 1655341 & 45/24.98/1 & 133.48 & & 37/65.30/1 & 265.49 \\
\bottomrule[1.5pt]
\end{tabular}
\caption{Number of iterations of P-MINRES, mean number of iterations for the PCG($\An$) and PCG($Q$), and total time for $k'=3$ (lid-driven cavity).} 
\label{tab:ic_pgc_tols_k3} 
\end{table}

As can be seen by the second column, $R+F(s)$, of Tables \ref{tab:ic_pgc_tols_k2} and \ref{tab:ic_pgc_tols_k3} that measures the time spent on the reordering and factorization steps these procedures worked extremely fast in both cases, when compared to the overall solver time. With respect to the number of non-zeros of the factors for $\An$ and $Q$, one can see that, for both cases and all mesh-sizes, they are almost half the number of non-zeros of $\An$ and $Q$ themselves. Again we see a practically constant number of iterations for P-MINRES for both $k'=2$ and $k'=3$, and here also for both $tol = 10^{-3}$ and $tol = 10^{-6}$, indicating a spectral equivalence with mesh refinement. 

Moreover, comparing the number of P-MINRES iterations of Tables \ref{tab:niter_cavity_k_2} and \ref{tab:neg_and_pos_eig_k_3}, column PCG($\An,Q$), with those for IC(0)-PCG($\An,Q$) in Tables \ref{tab:ic_pgc_tols_k2} and \ref{tab:ic_pgc_tols_k3}, we see they were almost the same, but obviously for the incomplete Cholesky case, that offers a better preconditioning, the number of internal iterations of PCG($\An$) and PCG($Q$) were considerably reduced. Surprisingly, PCG($Q$) converged in one iteration for all mesh-sizes and tolerances, showing that reordering by RCM and incomplete Cholesky offers a good preconditioning for the pressure mass matrix $Q$. Comparing the total time for PCG($\An,Q$) and IC(0)-PCG($\An,Q$) ($tol = 10^{-6}$) we observed a reduction of more than 10 seconds for $k'=2$, $h=1/128$ and almost 20 seconds for $k'=3$, $h=1/128$  

The relaxation of the relative residual tolerance for PCG also improved the total time for both $k'=2$ and $k'=3$. We can observe that the number of iterations of P-MINRES increased a little, on the other hand, the mean number of inner iterations of PCG($\An$) decreased, causing an overall decrease in time as compared to the case where $tol = 10^{-6}$.

Finally, we observe that the IC(0)-PCG($\An,Q$) preconditioning strategy gave the best time results of all strategies, losing only to the sparse direct solver when applied to $k'=2$, and the mesh-sizes $h=1/64$ and $h=1/128$. That is why we tested additionally the mesh-size $h=1/256$, for both IC(0)-PCG($\An,Q$) and the direct solver, where we see that IC(0)-PCG($\An,Q$) with PCG tolerance of $10^{-3}$ performed almost 4 times faster. Also, motivated by the excellent results of \cite{Kronbichler2012}, which more relaxed PCG tolerances were used, we ran case $k'=3$, $h=1/256$ with a PCG tolerance of $10^{-2}$ and the total time are $153.97$ seconds, which is bigger than our best time, 133.48 seconds with PCG tolerance of $10^{-3}$.

\section{Conclusion and future work}
Divergence-conforming B-spline discretizations are recent schemes based on the Isogeometric concept. In addition to being inf-sup stable, they are also divergence-free, a feature that is not easily achieved by mixed inf-sup stable elements, nor for stabilized ones. Their mathematical properties, presented by Evans \cite{Evans2011}, highlight their potential for viscous incompressible flows analyses. As usual, divergence-conforming discretizations end up in a linear system, and the efficient solution of such systems is of fundamental importance. In this paper, we analyze the performance of block-diagonal preconditioners, as introduced by Wathen and Silvester \cite{Wathen1993},\cite{Wathen1994} and \cite{elman2005finite}, for divergence-conforming discretizations, applied to the Stokes problem.

We have shown by means of several choices for the blocks that the theoretical bounds for the spectra, derived by Wathen and Silvester \cite{Wathen1993}, \cite{Wathen1994}, in the context of classical elements, also holds for divergence-conforming discretizations. Moreover, the linear system coefficient matrix spectrum allows the construction of spectrally equivalent preconditioners for Stokes problems. One of the ingredients, in the block-diagonal preconditioning strategy, is the proper approximation of the pressure mass matrix $Q$. We have shown that, for divergence-conforming discretizations, the approximation of $Q$ by taking its diagonal entries generates a poor preconditioner, particularly for the polynomial degree $k'=3$.

Another ingredient is the approximation of the viscosity matrix $\An$. As the lid-cavity flow results have shown, the use of iterative solvers to approximate $\An$, like preconditioned conjugate gradients, performed well for both polynomial degrees $k'=2$ and $k'=3$, given the smallest time results for fine meshes. Also, we have shown that reordering the unknowns with zero fill-in Incomplete Cholesky factorization as preconditioners, for both $\An$ and $Q$, with relaxed inner relative residual tolerances, yield very good preconditioners. Nevertheless, there is still room for improvements since the mean number of iterations of the inner preconditioned conjugate gradients applied to $\An$ is not spectrally equivalent with mesh refinement.

Several investigations may unfold. Since we only tested two-dimensional problems, we feel that the performance  evaluation of Krylov solvers and block preconditioning strategies for divergence-conforming discretizations applied to larger problems, followed by a scalability analysis, must be done. Another aspect that we like to pursue is the coupling of incompressible flows, discretized by divergence-conforming spaces, with transport problems.

\bibliographystyle{unsrt}
\bibliography{bib_thesis}

\end{document}